\documentclass[a4paper,twoside,12pt]{article}
\usepackage[english]{babel}

\usepackage{amsmath,amsfonts,amssymb,amsthm,graphicx}
\newtheorem{teo}{Theorem}
\newtheorem{lem}{Lemma}
\newtheorem{deff}{Definition}

 \title{{\bf  On the  essential decreasing   of the summation order in the Abel-Lidskii sense     }}

\author{Maksim \,V.~Kukushkin   \\ \\
  \small\textit{Russian Academy of Sciences,  Kabardino-Balkarian Scientific Center,}\\
 \small\textit{Institute of Applied Mathematics and Automation, 360000,  Nalchik, Russia}\\
 \small  \textit{National Research University Higher School of Economics, 101000,  Moscow, Russia}\\}

\date{}
\begin{document}
\maketitle

\begin{abstract} In this  paper, we consider a problem of decreasing the summation order in the Abel-Lidskii sense. The problem has a significant prehistory since 1962 created by such mathematicians as   Lidskii V.B.,
Katsnelson V.E., Matsaev V.I.,     Agranovich M.S. As a main result, we  will show that the summation order  can be decreased   from the values more than a convergence exponent, in accordance with the Lidskii V.B. results, to   an arbitrary small positive number. Additionally, we construct a qualitative theory of summation in the Abel-Lidkii sense and produce a number of fundamental  propositions that may represent the interest themselves.

\end{abstract}
\begin{small}\textbf{Keywords:}
   Abel-Lidskii basis property;   Schatten-von Neumann  class; convergence exponent; counting function; sectorial operator.   \\\\
{\textbf{MSC} 47B28; 46B15; 47A10; 47B12; 47B10;   58D25.}
\end{small}

\section{Introduction}

\subsection{Historical review}

$\quad$   In order to establish  a harmonious connection between   well-known  facts,  recall   that the eigenvectors system of a compact selfadjoint operator forms a basis in the closure of its range. This fact can be interpreted in terms of the  spectral theorem as a statement on the unit decomposition  \cite{firstab_lit:sadovn}.  Consider    a more  general case corresponding to a compact  non-selfadjoint operator   with  the numerical range of values   belonging  to a sector with the semi-angle less than $\pi$  and the vertex situated at the point zero \cite{firstab_lit:1Lidskii}. Obviously, the case covers a compact non-negative selfadjoint operator. Apparently,   we cannot weaken conditions upon the numerical range   expecting  that the basis property   would be preserved, moreover the fact of  the root vectors system  completeness becomes non-obvious \cite{firstab_lit:1Krein},\cite{firstab_lit:2Lidskii} what makes   a prerequisite for a comprehensive study of the issue.

In the recent century,   the problem on  root vectors system completeness  related to   non-selfadjoint operators  attracted  a serious attention. The beginning of the research  was laid in the paper by  Keldysh M.V. 1959 \cite{firstab_lit:Keldysh M.V.}. In subsequent works, there were found sufficient conditions for completeness of the root vectors system 1958-1959 \cite{firstab_lit:1958Lidskii},\cite{firstab_lit:1959Lidskii}, \cite{firstab_lit:1959Krein},\cite{firstab_lit:1959Sakhnovich}. However, the fact is that   the  completeness property   of the root vectors system   is not sufficient for  the basis property.   This fact establishes  a prerequisite for creation  of  a  fundamental direction  in  the abstract spectral theory devoted to the basis property of the root vector system in a generalized sense. In the papers by Markus A.S. 1960 \cite{firstab_lit:1Markus}, 1962 \cite{firstab_lit:2Markus} the problem on the series convergence in the Bari sense of subspaces was considered. In the papers  1960 \cite{firstab_lit:1960Lidskii}, 1962 \cite{firstab_lit:1Lidskii}   Lidskii V.B. introduced a generalized $(A,\lambda,s)$ - method  of summation  for series on the root vectors based on the notion of the Abelian means considered in the monograph by  Hardy G.H. 1949 \cite{firstab_lit:Hardy},  the parameter $s$ is called by the order of summation. The generalization for Banach spaces was considered by Markus A.S. in 1966 \cite{firstab_lit:1966Markus}.  In the paper by Agranovich M.S. 1976 \cite{firstab_lit:Agranovich1976} a class of non-selfadjoint elliptic operators was considered in the framework of the problem. The formula connecting spectral asymptotics corresponding to an operator with a discrete spectrum  and its real component was established by  Markus A.S.,  Matsaev V.I. in 1981 \cite{firstab_lit:Markus Matsaev}, the authors established the sufficient conditions for unconditional basis property in the sense of subspaces 1981 \cite{firstab_lit:Markus Matsaev1}. The problem on preservation of the unconditional basis property under non-selfadjoint perturbations of selfadjoint operators was considered by Motovilov A.K., Shkalikov A.A. in the paper 2019 \cite{firstab_lit:Motovilov}. The latest overview of the results related to the problem on decomposition on the  root vectors series   was  represented  in the paper by Shkalikov A.A. 2016 \cite{firstab_lit:Shkalikov A.}.

The problem on decreasing of the summation order   was firstly formulated   by Lidskii V.B. in the paper  1962 \cite{firstab_lit:3Lidskii}  not for the general case but for the case corresponding to the perturbation of the selfadjoint elliptic operator of the second order under  the strong  subordination condition. More generally,  the problem was considered by  Katsnelson V.E. in the Ph.D thesis 1967 \cite{firstab_lit:1967Katsnelson} (see also \cite{firstab_lit: 1999 Katsenelenbaum}) for  perturbations of non-negative selfadjoint operators. The problem on  decreasing of the summation order for  operators with the spectrum belonging to the domain of the parabolic and hyperbolic  type was considered   by Shkalikov A.A. in 1982 \cite{firstab_lit:1982 Shkalikov}, 1983 \cite{firstab_lit:1983 Shkalikov}. The detailed   substantiation of a method allowing to decrease the summation order was represented by Agranovich M.S. in 1994 \cite{firstab_lit:2Agranovich1994} for operators with the numerical range of values belonging to a domain of the parabolic type. Moreover,   a general scale of conditions admitting    convergence of the root vector series in a generalized sense such as  Bari, Riesz, Abel-Lidskii senses of the series convergence was established. The clarification of the summation order  was  implemented    in the paper  by Kukushkin M.V. 2022 \cite{firstab_lit:1kukushkin2021}.

Apparently, the main advantage of the  Lidskii V.B. \cite{firstab_lit:1Lidskii}  method    is  wider assumptions related to the numerical range of values  comparatively with the sectorial condition let alone the operator class corresponding  the numerical range of values  belonging to the domain of the parabolic type \cite{firstab_lit:2Agranovich1994},\cite{firstab_lit:Markus Matsaev},\cite{firstab_lit:2Markus}. Note that such a location of the numerical range of values is the inherent property of the operators with a salfadjoint senior term. At the same time, a scientific novelty and  relevance appear     in  the very case when  a senior term  is not selfadjoint  for there exists a comprehensive theory  devoted to   perturbed selfadjoint operators, see papers
\cite{firstab_lit:1Katsnelson},\cite{firstab_lit:1Krein},\cite{firstab_lit:Markus Matsaev},\cite{firstab_lit:1Markus},\cite{firstab_lit:2Markus},\cite{firstab_lit:Motovilov},\cite{firstab_lit:Shkalikov A.}.  The fact is that most of them deal with a decomposition of the  operator  on a sum,  where the senior term  must be either a selfadjoint or normal operator. Otherwise, the  methods of the papers
     \cite{kukushkin2019,firstab_lit(arXiv non-self)kukushkin2018} become relevant   since they  allow us to study spectral properties of   operators  whether we have the mentioned above  representation or not, moreover   they have a natural mathematical origin that appears brightly  while we are considering abstract constructions expressed  in terms of the semigroup theory   \cite{kukushkin2021a}.

  In this paper we consider a sectorial operator belonging to the trace class. Generally, we will show that the summation order  can be decreased to the value depending on the growth of the  algebraic multiplicities. In particular, we establish a remarkable  fact that the summation order is an arbitrary small number for operators having a sufficiently low growth of the algebraic multiplicities. The  idea to consider the growth of the algebraic multiplicities as an inherent  property is principally novel in comparison with the previously obtained results and  the fact that the summation order can be decreased for a sectorial operator is  new.        The main application of this paper  results appeals to the abstract  Cauchy problem for the  evolution equation  including the qualitative theory for   fractional evolution equations   \cite{firstab_lit:Bazhl},\cite{firstab_lit:Bazhl1},\cite{firstab_lit:Bazhl2}. More detailed survey related to applications to various  physical-chemical processes  and   applications to well-known concrete  evolution equations of the fractional order is given in the papers   \cite{kukushkin2019},\cite{kukushkin2021a},\cite{firstab_lit:1kukushkin2021},\cite{firstab_lit:2kukushkin2022},\cite{firstab_lit(axi2022)}.   Undoubtedly, the  main achievement  of this paper is a  constructed abstract qualitative theory creating an opportunity to solve more concrete problems let alone
far reaching modifications and  generalizations.

\subsection{Preliminaries}

Let    $ C,C_{i} ,\;i\in \mathbb{N}_{0}$ be   real positive  constants. We   assume   that     values of $C$    can be different in     formulas  but   values of $C_{i} $ are  certain.    Everywhere further,   we consider   linear    densely defined operators acting in a separable complex  Hilbert space $\mathfrak{H}$. Denote by $ \mathcal{B} (\mathfrak{H})$    the set of linear bounded operators   on    $\mathfrak{H}.$     Denote by    $    \mathrm{D}   (L),\,   \mathrm{R}   (L),\,\mathrm{N}(L),\,\mathrm{P}(L)$   the  {\it domain of definition}, the {\it range},  the {\it kernel}, and   the {\it resolvent set}  of the  operator $L$ respectively. Denote by $\Sigma(L):=\mathbb{C}\backslash \mathrm{P}(L)$ the spectrum of the operator $L.$

Consider a pair of complex Hilbert spaces $\mathfrak{H},\mathfrak{H}_{+},$ the notation
$
\mathfrak{H}_{+}\subset\subset\mathfrak{ H}
$
means that $\mathfrak{H}_{+}$ is dense in $\mathfrak{H}$ as a set of    elements and we have a bounded embedding provided by the inequality
$$
\|f\|_{\mathfrak{H}}\leq C\|f\|_{\mathfrak{H}_{+}},\;f\in \mathfrak{H}_{+},
$$
moreover   any  bounded  set with respect to the norm $\mathfrak{H}_{+}$ is compact with respect to the norm $\mathfrak{H}.$
Denote by  $ \mathfrak{Re} L  := \left(L+L^{*}\right)/2,\, \mathfrak{Im} L  := \left(L-L^{*}\right)/2 i$
  the  real  and   imaginary components    of the  operator $L$  respectively.
In accordance with  the terminology of the monograph  \cite{firstab_lit:kato1980} the set $\Theta(L):=\{z\in \mathbb{C}: z=(Lf,f)_{\mathfrak{H}},\,f\in  \mathrm{D} (L),\,\|f\|_{\mathfrak{H}}=1\}$ is called the  {\it numerical range}  of the   operator $L.$
Define a closed sector in the complex  plain   $\mathfrak{ L}_{a}(\theta):=\{z\in \mathbb{C}:  |\arg(z-a)|\leq\theta<\pi\}\cup \{a\},$ where $a\in \mathbb{C}$ is called by  the vertex and $\theta$ is called by  the semi-angle of the sector.     An  operator $L$ is said to be     {\it sectorial}    if its  numerical range   belongs to a
sector     $\mathfrak{ L}_{a}(\theta),\;\theta<\pi/2.$
   An operator $L$ is said to be   {\it   accretive}   if  $\mathrm{Re}(Lf,f)_{\mathfrak{H}}\geq 0,\,f\in \mathrm{D}(L),$  {\it m-accretive} if  $\{z\in \mathbb{C}:\,\mathrm{Re} z<0\}\subset \mathrm{P}(L),\,\|(L-\lambda I)^{-1}\|\leq  |\mathrm{Re} \lambda|^{-1},\,\mathrm{Re }\lambda<0,$ {\it   dissipative}
   if  $\mathrm{Im}(Lf,f)_{\mathfrak{H}}\geq 0,\,f\in \mathrm{D}(L).$
    An operator $L$ is called by the  {\it operator with  discrete spectrum} if  $0\in \mathrm{P}(L)$ and the inverse operator is compact.
The  dimension of the root vectors subspace corresponding to a certain  eigenvalue  of the operator $L$ is called by the {\it algebraic multiplicity} of the eigenvalue. Denote by $\mu_{j}(L),\,j\in \mathbb{N} $ the eigenvalues of the operator $L,$ where the numbering  is given in accordance with   increase  or   decrease of their absolute value and with this numbering each eigenvalue is counted as many times as its algebraic multiplicity.   Denote by   $\nu(\mu_{j})$ the  algebraic multiplicity of the eigenvalue $\mu_{j}(L)$  and denote  by  $\nu(L)$  the sum of all algebraic multiplicities of the  operator $L.$
     Suppose $L$ is  a compact operator and  $|L|:=(L^{\ast}L)^{1/2},$ then   the eigenvalues of the operator $|L|$ are called by  the {\it singular  numbers} ({\it s-numbers}) of the operator $L$ and are denoted by $s_{j}(L),\,j=1,\,2,...\,,{\rm dim}\,  \mathrm{R}  (|L|).$ If ${\rm dim}\,  \mathrm{R}  (|L|)<\infty,$ then we put by definition     $s_{j}=0,\,j>{\rm dim}\,  \mathrm{R}  (|L|).$
 Assume that an  operator $L$ is compact, the following relation holds
$$
\sum\limits_{n=1}^{\infty}s^{\sigma}_{n}(L)<\infty,\,0<\sigma<\infty,
$$
then $L$ is said to be in the Schatten-von Neumann  class $\mathfrak{S}_{\sigma}(\mathfrak{H}),$ i.e. $L\in \mathfrak{S}_{\sigma}(\mathfrak{H})$ in symbol.  Denote by   $\mathfrak{S}_{\infty}(\mathfrak{H})$ the  set of compact operators acting in $\mathfrak{H}.$
Let $L$ be a bounded operator acting in $\mathfrak{H},$ and assume that $\{\varphi_{n}\}_{1}^{\infty},\,\{\psi_{n}\}_{1}^{\infty}$ a pair of orthonormal bases in $\mathfrak{H}.$ Define the {\it absolute operator norm} as follows
$$
\|L\|_{2}:=\left(\sum\limits_{n,k=1}^{\infty}|(L\varphi_{n},\psi_{k})_{\mathfrak{H}}|^{2}\right)^{1/2}.
$$
 Consider a sequence   $\{a_{j}\}_{1}^{\infty}\subset \mathbb{C},$ define  the    counting function
$
n(r,a_{j}):=\mathrm{card}\{j\in \mathbb{N}:|a_{j}|\leq r\}.
$
Assume that an operator  $L$ is compact (operator with discrete spectrum), denote by $n(r,L)$  the counting function corresponding to the sequence of the absolute values of the  operator  characteristic numbers (eigenvalues).

 Further, we consider   a compact    operator $B,$   observe  the sequence of its eigenvalues
$$
0\neq\mu_{1}=...=\mu_{p_{1}}\neq \mu_{p_{1}+1}=...=\mu_{p_{2}}\neq \mu_{p_{2}+1}=...=\mu_{p_{3}}\neq...\,.
$$
  Analogously to the definitions accepted in the entire function theory \cite{firstab_lit:Baraichev}, we will call the numbers  $p_{j},\,j\in \mathbb{N}$  by   the {\it  principal indexes}. In accordance with the above, we have
$$
\Delta_{j}= p_{j} - p_{j-1}
,\,j\in \mathbb{N},
$$
where    $\Delta_{j}$ denotes the algebraic multiplicity  of the  eigenvalue     $\mu_{p_{j}},$ we formally  put   $p_{0}=0.$ Denote by $\lambda_{j}: =1/\mu_{j},\,j\in \mathbb{N}$ the  {\it characteristic numbers} of the operator $B.$ For the reader convenience,   we use    special   notations  for the  eigenvalues and the characteristic numbers corresponding to the principal indexes  $y_{j}:=\mu_{p_{j}},\,z_{j}:= \lambda_{p_{j}},\,j\in \mathbb{N},$  we  will call them by the  {\it principal eigenvalues} and   the {\it principal characteristic  numbers}  of the operator $B$ respectively.
    We also follow the definitions and  notations accepted in  the monographs    \cite{firstab_lit:1Gohberg1965},  \cite{firstab_lit:kato1980}.\\

\section{Overview of the supplementary results}

    Following  the monograph \cite{firstab_lit:Eb. Levin}, we introduce some notions and facts of the entire function theory. In this subsection, we   use the following notations
$$
G(z,p):=(1-z)e^{z+\frac{z^{2}}{2}+...+\frac{z^{p}}{p}},\,G(z,0):=(1-z).
$$
Consider   an entire function  that has zeros satisfying the following relation for some   $\beta>0$
\begin{equation}\label{1}
 \sum\limits_{n=1}^{\infty}\frac{1}{|a_{n}|^{\beta}}<\infty.
\end{equation}
In this case, we denote by $p$ the smallest integer number for which the following condition holds
\begin{equation}\label{2}
\sum\limits_{n=1}^{\infty}\frac{1}{|a_{n}|^{p+1}}<\infty .
\end{equation}
 It is clear that $0\leq p<\beta.$ Consider a formal infinite product
 \begin{equation}\label{3}
  \prod\limits_{n=1}^{\infty} G\left(\frac{z}{a_{n}},p\right),
\end{equation}
   we will call it a canonical product and call $p$ the genus of the canonical product.
By the   {\it convergence exponent} $\rho$ of the sequence
$
\{a_{n}\}_{1}^{\infty}\subset \mathbb{C},\,a_{n}\neq 0,\,a_{n}\rightarrow \infty
$
 we mean the greatest lower bound for such numbers $\beta$ that the   series \eqref{1} converges.
   We need the following lemma, see Theorem 2 \cite[p.29]{firstab_lit:Eb. LevinE}, Lemma 3 \cite{firstab_lit:Eb. Levin}.
\begin{lem} \label{L1}
 If the series \eqref{2} converges, then the corresponding infinite product \eqref{3} converges uniformly on every compact subset and satisfies  the estimate
 $$
 \ln\left|  \prod\limits_{n=1}^{\infty} G\left(\frac{z}{a_{n}},p\right)\right|\leq k_{p}r^{p}\left(\int\limits_{0}^{r}\frac{n(t)}{t^{p+1}}dt+r\int\limits_{r}^{\infty}\frac{n(t)}{t^{p+2}}dt\right),\,r:=|z|,
 $$
 where $k_{p}=3e(p+1)(2+\ln p),\,p>0,\,k_{0}=1.$
\end{lem}

\subsection{Characteristic  determinant}\label{p 2.1}

The well-known technique used by Lidskii V.B. \cite{firstab_lit:1Lidskii} and others appeals to the notion of the characteristic  determinant and due to this reason  we produce a complete description of the object. Having chosen an orthonormal basis $\{e_{j}\}_{1}^{\infty}\subset \mathfrak{H}$ consider a matrix $\{b_{ij}\}_{1}^{\infty}$    of the operator $B\in \mathfrak{S}_{\infty},$  where
$$
b_{ij}:=(Be_{j,}e_{i})_{\mathfrak{H}},\,i,j=1,2,...\,.
$$
Assume that the finite-dimensional  space $\mathbb{E}^{n}$ generated by the vectors $\{e_{j}\}_{1}^{n}$ is an invariant space of the operator $B,$ thus we have a restriction $B_{n}\subset B,\,B_{n}:\mathbb{E}^{n}\rightarrow \mathbb{E}^{n}.$  Denote by $\det\{I-B_{n}\}$ the determinant of the matrix $\{\delta_{ij}-  b_{ij}\}_{1}^{n}.$ It is a well-known fact that   the determinant $\det\{I-B_{n}\}$ does not depend on a basis in $\mathbb{E}^{n}$ since
$$
\det\{I-B_{n}\}=\prod\limits_{j=1}^{\nu(B_{n})}\left(1-\mu_{j}(B_{n})\right),
$$
where $\nu(B_{n})$ is the algebraic multiplicity (dimension of the root vector subspace) of the operator $B_{n}.$ The latter relation  shows that it is possible to make sense for the following construction
$$
\det\{I-B\}=\prod\limits_{j=1}^{\nu(B)}\left(1-\mu_{j}(B)\right),\, B\in \mathfrak{S}_{1},
$$
where $\nu(B)\leq\infty.$ The product at the right-hand side of the last relation is convergent since, we have
$$
 B\in \mathfrak{S}_{1} \Rightarrow \sum\limits_{j=1}^{\nu(B)}|\mu_{j}(B)|<\infty.
$$
 Now, consider a formal decomposition of the determinant of the matrix with the infinite quantity of rows and columns
$$
\Delta(\lambda):=\det \{\delta_{ij}-\lambda b_{ij}\}_{1}^{\infty}=\sum\limits_{p=0}^{\infty}(-1)^{p}q_{p}\lambda^{p},\,\lambda\in \mathbb{C},
$$
where $q_{0}=1$ and $q_{p},\,p=1,2,...\,$ is a sum of all central minors of the matrix $\{b_{ij}\}_{1}^{\infty}$  of the order   $p,$   formed from the columns and rows with $i_{1},i_{2},...,i_{p}$ numbers,  i.e.
$$
q_{p}=\frac{1}{p!}\sum\limits_{i_{1},i_{2},...,i_{p} =1}^{\infty}B\begin{pmatrix} i_{1}&  i_{2}&...&i_{p}\\
  i_{1}&  i_{2}&...&i_{p}
\end{pmatrix}.
$$
Note that if $B\in \mathfrak{S}_{1}$   then in accordance with the well-known theorems (see \cite{firstab_lit:1Gohberg1965}), we have
\begin{equation}\label{4}
\sum\limits_{n  =1}^{\infty}|b_{nn}| <\infty,\;\sum\limits_{n,m  =1}^{\infty}|b_{nm}|^{2} <\infty,
\end{equation}
where $ b_{nm} $ is the matrix coefficients of the operator $B.$ This follows easily from the properties of the trace class operators and Hilbert-Schmidtt class operators respectively. In accordance with the von Koch H. theorem \cite{firstab_lit:Koch1909}  conditions \eqref{4} guaranty the absolute convergence of the series $q_{p}.$
Moreover, the formal series $\Delta(\lambda)$ is convergent for arbitrary $\lambda\in \mathbb{C},$  therefore  it represents an entire function. Analogous facts take place if we consider a formal decomposition of a minor  corresponding to the matrix obtained due to deleting   the  $l$ -th row and the $m$ -th column from the initial matrix $\{\delta_{ij}-\lambda b_{ij}\}_{1}^{\infty}.$ Thus, we can give  a  meaning  to the following construction
$$
\Delta^{lm}(\lambda):= 1+\sum\limits_{p=1}^{\infty}(-1)^{p}\lambda^{p}\!\!\!\!\!\sum\limits_{i_{1},i_{2},...,i_{p} =1}^{\infty}\!\!\!B\begin{pmatrix} i_{1}&  i_{2}&...&i_{p}\\
  i_{1}&  i_{2}&...&i_{p}
\end{pmatrix}_{lm},\,\lambda\in \mathbb{C},
$$
where the  used formula in brackets means a minor    formed from the columns and rows with $i_{1},i_{2},...,i_{p}$ numbers corresponding to the matrix obtained due to deleting   the  $l$-th row and the $m$-th column from the initial matrix.  Conditions \eqref{4} guarantee   convergence of  formal series $\Delta^{lm}(\lambda)$  for an arbitrary $\lambda\in \mathbb{C}.$

Assume that $\lambda$ is a regular point of the operator  $(I-\lambda B)^{-1}$  in accordance with (1.11) \cite{firstab_lit:1Lidskii}   the  equation
$
(I-\lambda B)x=f,
$
where $x,f\in \mathfrak{H},$
can be rewritten as a system in the   form
\begin{equation*}
\sum\limits_{j=1}^{\infty}(\delta_{ij}-\lambda b_{ij})x_{j}=f_{i},\;f_{i}=(f,e_{i})_{\mathfrak{H}},\,i=1,2,...\,.
\end{equation*}
In accordance with \cite{firstab_lit:1Lidskii}   conditions \eqref{4} guarantee existence of the solution in the form
\begin{equation}\label{5}
(I-\lambda B)^{-1}f =\sum\limits_{m=1}^{\infty}\left(\sum\limits_{l=1}^{\infty}(-1)^{l+m}\frac{\Delta^{lm}(\lambda)}{\Delta(\lambda)}f_{l}\right)e_{m},
\end{equation}
where $f_{l}=(f,e_{l})_{\mathfrak{H}}.$ The entire function $\Delta(\lambda)$ is called  by  the   {\it Fredholm determinant} of the operator $B.$ In accordance with the definition \S 4, Chapter I, \cite{firstab_lit:1Gohberg1965} under assumption $B\in \mathfrak{S}_{1}$ the product  $\det\{I-\lambda B\}$ is called by  the {\it characteristic determinant} of the operator $B.$

Since the main characteristic of the studied operators is the   Schatten   classification  then it is rather reasonable to  provide  auxiliary   propositions formulated in corresponding   terms. The following lemma is represented in  \cite{firstab_lit:1Lidskii}.

\begin{lem}\label{L2}
  Assume that $B$ is a compact operator, $P$ is an arbitrary orthogonal projector in $\mathfrak{H},$  then
$$
s_{n}(PBP)\leq s_{n}(B),\,n\in \mathbb{N}.
$$
\end{lem}

The statement  of the  following lemma is included  in  the proof of Lemma 2 \cite{firstab_lit:1Lidskii}, here for the reader convenience, we represent it  supplied with   expended reasonings.

\begin{lem}\label{L3} 
Assume that  $B\in \mathfrak{S}_{1},$  then the following representation holds
$$
\Delta(\lambda)=\prod\limits_{n=1}^{\infty}\left\{1- \lambda \mu_{n}(B)    \right\},\,\lambda\in \mathbb{C},\,i.e. \;\Delta(\lambda)=\det\{I-\lambda B\},
$$
the characteristic and Fredholm determinants of the operator $B$ are coincided.
\end{lem}
\begin{proof} Firstly, we should note that in accordance with Theorem 8.1 \S 8, Chapter III, \cite{firstab_lit:1Gohberg1965} the operator $B$  belongs to the trace class. Therefore,  for an arbitrary orthonormal basis $\{\varphi_{n}\}_{1}^{\infty},$  we have
$$
\sum\limits_{n=1}^{\infty}(B\varphi_{n},\varphi_{n})_{\mathfrak{H}}<\infty.
$$
  The arbitrariness in the choice
of a basis gives an opportunity to claim  that the series is convergent after an arbitrary transposition of the terms, from what follows that the series is absolutely convergent. Hence, the first condition \eqref{4} holds. To prove the second condition \eqref{4}, we should note   the inclusion $\mathfrak{S}_{2}\subset \mathfrak{S}_{1}$ and the fact that $\mathfrak{S}_{2}$ coincides with the so-called Schmidt class of  operators having a finite absolute operator norm $\|\cdot\|_{2}.$ The latter fact can   be established if we consider a complement  of the orthonormal set $\{\varphi_{n}\}_{1}^{\infty}$ of the eigenvectors of the operator $B^{\ast}B$ to a basis $\{\psi_{n}\}_{1}^{\infty}$ in the Hilbert space. Then in accordance with the well-known decomposition formula (see \S 3, Chapter V,  \cite{firstab_lit:kato1980}), we get the orthogonal sum
$$
\mathfrak{H}=\overline{\mathrm{R}(B^{\ast}B)}\,\dot{+}\,\mathrm{N}(B^{\ast}B),
$$
where $\{\varphi_{n}\}_{1}^{\infty}$ is a basis in  $\overline{\mathrm{R}(B^{\ast}B)},$ in accordance with the well-known  property of compact selfadjoint operators. Therefore, a complement  of the system $\{\varphi_{n}\}_{1}^{\infty}$ to the basis in $\mathfrak{H}$ belongs to $\mathrm{N}(B^{\ast}B).$ Hence
$$
\|B\|^{2}_{2}=\sum\limits_{n,k=1}^{\infty}|(B\psi_{n},\psi_{k})_{\mathfrak{H}}|^{2} =\sum\limits_{n=1}^{\infty}\|B\psi_{n}\|^{2}_{\mathfrak{H}}=\sum\limits_{n=1}^{\infty}(B^{\ast}B\psi_{n}, \psi_{n})_{\mathfrak{H}}  =\sum\limits_{n=1}^{\infty}(B^{\ast}B\varphi_{n}, \varphi_{n})_{\mathfrak{H}}  =\sum\limits_{n=1}^{\infty}s^{2}_{n}.
$$
 Therefore, by virtue of belonging to the trace class    conditions \eqref{4} hold for an arbitrary chosen basis in the Hilbert space.

  Now, consider an arbitrary basis $\{\varphi_{k}\}^{\infty}_{1}$ and consider an orthogonal  projector $P_{n}$ corresponding to the  subspace generated by the first $n$ basis vectors $\varphi_{1},\varphi_{2},...,\varphi_{n}.$ Consider a determinant
  $$
  \Delta^{(n)}(\lambda):=\det \{\delta_{ij}-\lambda b_{ij}\}_{ij=1}^{n} =\prod\limits_{k=1}^{n}\left\{1-  \lambda \mu_{k}( P_{n}BP_{n} )  \right\}
  $$
Using the Weil inequalities \cite{firstab_lit:1Gohberg1965}, we get
$$
  |\Delta^{(n)}(\lambda)| \leq\prod\limits_{k=1}^{n}\left\{1+  |\lambda|\cdot |\mu_{k}( P_{n}BP_{n} )|  \right\}\leq \prod\limits_{k=1}^{n}\left\{1+  |\lambda|\cdot |s_{k}( P_{n}BP_{n} )|  \right\}.
  $$
Applying Lemma \ref{L2}, we get
$$
  |\Delta^{(n)}(\lambda)|  \leq \prod\limits_{k=1}^{\infty}\left\{1+  |\lambda|\cdot |s_{k}(  B  )|  \right\}.
  $$
Passing to the limit while $n\rightarrow\infty,$ we get
$$
  |\Delta(\lambda)|  \leq \prod\limits_{k=1}^{\infty}\left\{1+  |\lambda|\cdot |s_{k}(  B  )|  \right\}.
  $$
It implies, if we observe Theorem 4 (Chapter I, \S 4) \cite{firstab_lit:Eb. Levin} that the entire function $\Delta (\lambda)$ is of the finite order. Therefore, in accordance with  Theorem 13 (Chapter I, \S 10) \cite{firstab_lit:Eb. Levin}, it has a representation by the canonical product, i.e.
\begin{equation*}
\Delta (\lambda)= \prod\limits_{n=1}^{\infty}\left\{1- \lambda \mu_{n}(B )    \right\}.
\end{equation*}
The proof is complete.
\end{proof}

\subsection{ Abel-Lidskii   Series expansion}\label{P2.2}

   In accordance with the Hilbert theorem
  (see \cite{firstab_lit:Riesz1955}, \cite[p.32]{firstab_lit:1Gohberg1965})   the spectrum of an arbitrary  compact operator $B$  consists of the so-called normal eigenvalues, it gives us an opportunity to consider a decomposition to a direct sum of subspaces
\begin{equation}\label{6}
 \mathfrak{H}=\mathfrak{N}_{q} \oplus  \mathfrak{M}_{q},
 \end{equation}
corresponding to the principal eigenvalue   $y_{q},\,q\in \mathbb{N},$
 both  summands are   invariant subspaces of the operator $B,$  the first one is  a finite dimensional root subspace corresponding to the eigenvalue
 $y_{q},\, \mathrm{dim}\, \mathfrak{N}_{q}=\nu(y_{q}) $ and the second one is a subspace  wherein the operator  $B-y_{q} I$ is invertible.     We can choose the  Jordan basis  in $\mathfrak{N}_{q}$ that consists of Jordan chains of eigenvectors and root vectors  of the operator $B.$ Considering the set  of the principal   eigenvalues    and  corresponding    Jordan bases, we can arrange  a  {\it root vectors system}  or following  to the definition by   Lidskii V.B. \cite{firstab_lit:1Lidskii}  a system of the {\it major vectors} of the operator $B.$ The Riesz integral operator is  defined as follows
$$
\mathcal{P}_{ q }f:= -\frac{1}{2\pi i}\oint_{\Gamma'_{q}}(B-\lambda I)^{-1}f\,d\lambda,\; f\in \mathfrak{H},
$$
where $\Gamma'_{q}$ is a closed contour bounding a domain containing the eigenvalue  $y_{q}$ only. The properties of the Riesz integral operator   are described in detail in  \S 1.3, Chapter I,   \cite{firstab_lit:1Gohberg1965}.
Consider the integral operator introduced by Lidskii V.B.
$$
\mathcal{P} _{q}(s,t)f=-\oint\limits_{\Gamma_{q}}   e^{- \lambda^{s}  t}B \left(I-\lambda B \right)^{-1}f d\lambda,\,f\in \mathfrak{H},\,s,t>0,
$$ where $\Gamma_{q}$ is a closed contour bounding a domain containing the characteristic number   $z_{q}$ only.
  The following  fact was established by Lidskii V.B. in Lemma 5 \cite{firstab_lit:1Lidskii}.
\begin{lem}\label{L3a} Assume that  $B\in \mathfrak{S}_{\infty},\, \Theta(B)\subset \mathfrak{L}_{0}(\theta),\,\theta<\pi,$  then
\begin{equation}\label{9}
 \frac{1}{2\pi i}\int\limits_{\vartheta} e^{- \lambda^{s}  t} B\left(I-\lambda B\right)^{-1}f d\lambda \rightarrow f,\,t\rightarrow 0,\;f\in \mathrm{R}(B),\,s>0,
\end{equation}
\end{lem}
We should stress that the proof  represented by Lidskii V.B. is true for an arbitrary small positive  value $s.$
 The proof corresponding to the case $f\in \mathfrak{H},$ under the special condition on the norm of the resolvent, is represented in \cite{firstab_lit: 1999 Katsenelenbaum},   the idea of the proof    can be found  in  Theorem 5.1 \cite{firstab_lit:Shkalikov A.}.
Below, we provide a detailed proof in terms of the operator with a discrete spectrum.

\begin{lem}\label{L5} Assume  that  $W$ is an operator   with   discrete spectrum $\Theta(W)\subset \mathfrak{L}_{0}(\theta),\,\theta<\pi,$  then the following relation holds
$$
f(t):= \frac{1}{2\pi i}\int\limits_{\vartheta}  e^{- \lambda^{s}  t} \left(W-\lambda I\right)^{-1}f d\lambda \stackrel{\mathfrak{H}\;}{\rightarrow} f,\,t\rightarrow 0,\;f\in \mathfrak{H},\,0<s<\pi/2\theta.
$$
\end{lem}

\begin{proof}
 Let us prove that
 $$
 \|(W-\lambda I)^{-1}\|\leq C|\lambda|^{-1},\,\lambda\in\{z\in \mathbb{C}:\,  \arg  z =\psi\},\,\theta<|\psi|<\pi/2.
 $$
 The inequality for the resolvent  holds by virtue of    Theorem 3.2   \cite[p.268]{firstab_lit:kato1980}, since $\Theta(W)\subset \mathfrak{L}_{0}(\theta),$ and as a result
$$
 \|(W-\lambda I)^{-1}\|\leq  (\mathrm{dist}   \{\lambda,\overline{\Theta(W)} \})^{-1} \leq (\mathrm{dist}\left\{\lambda,\mathfrak{L}_{0}(\theta)\right\})^{-1}
 =  \{|\lambda|  \sin(|\psi|-\theta)\}^{-1}    ,
 $$
 $$
 \lambda\in\{z\in \mathbb{C}:\,\mathrm{arg} z=\psi\}.
$$
Note that $\mathrm{D}(W)$ is dense  in $\mathfrak{H},$ therefore for an arbitrary element $f\in \mathfrak{H},$ we can choose a sequence $\{f_{n}\}_{1}^{\infty}\subset\mathrm{D}(W)$ such that
$$
f_{n}\stackrel{\mathfrak{H}\;}{\rightarrow}f.
$$
In accordance with  \eqref{9}, we have
$$
 f_{n}(t):=\frac{1}{2\pi i}\int\limits_{\vartheta  }  e^{- \lambda^{s}  t}(W-\lambda I)^{-1}f_{n} d\lambda \stackrel{\mathfrak{H}\;}{\rightarrow} f_{n},\,t\rightarrow 0.
$$
Consider the inequality 
$$
\|f(t)-f\|_{\mathfrak{H}}\leq \|f_{n}(t)-f_{n}\|_{\mathfrak{H}}+\|f_{n}(t)-f(t)\|_{\mathfrak{H}}+\|f_{n}-f\|_{\mathfrak{H}},\,f\in \mathfrak{H}.
$$
Thus, if we show that
\begin{equation}\label{10}
f_{n}(t)\rightrightarrows f(t),\,n\rightarrow \infty,
\end{equation}
i.e. the sequence  converges uniformly with respect to $t,$ then  we obtain the desired result.
Let us make a change of the variable $\lambda=\xi t^{-1/s},$ then the contour $\vartheta$ has  undergone to a transformation leading to a contour $\vartheta'$   with the same orientation and preserved tendency to the infinitely-distant point, we have
$$
2\pi\| f_{n}(t)-f(t)\|_{\mathfrak{H}}=  t^{-1/s}\left\|\int\limits_{\vartheta'  }  e^{- \xi^{s}   }(W-  t^{-1/s}\xi I)^{-1}(f_{n}-f) d\xi\right\|_{\mathfrak{H}}\leq
$$
$$
\leq t^{-1/s}\|f_{n}-f\|_{\mathfrak{H}} \int\limits_{\vartheta'  }  e^{- \mathrm{Re}\xi^{s}   }\left\|(W-  t^{-1/s}\xi I)^{-1} \right\|  |d\xi|
\leq C t^{-1/s}\|f_{n}-f\|_{\mathfrak{H}} \int\limits_{\vartheta'  }  e^{- \mathrm{Re}\xi^{s}   } t^{1/s} |\xi|^{-1} |d\xi|=I_{1}.
$$
Using the condition  $0<s<\pi/2\theta,$ we have  $\mathrm{Re}\xi^{s}>C|\xi|^{s},\,\xi\in \vartheta',$ therefore
$$
I_{1}\leq    \|f_{n}-f\|_{\mathfrak{H}} \int\limits_{\vartheta'  }  e^{- C |\xi|^{s}   }   |\xi|^{-1} |d\xi|\leq C\|f_{n}-f\|_{\mathfrak{H}}.
$$
The latter relation shows that \eqref{10} holds. The proof is complete.
\end{proof}

It is remarkable that the method for summation of the root vectors series    invented  by Lidskii V.B. \cite{firstab_lit:1Lidskii}  originates from the notion of the Abelian means considered in the monograph by Hardy G.H.  \cite[p.71]{firstab_lit:Hardy}.  We can apply  the original  definition in the following way. Consider a formal decomposition of an element $f\in \mathfrak{H}$  on the series
\begin{equation}\label{10y}
f\sim \sum\limits_{q=1}^{\infty}\mathcal{P}_{q}f.
\end{equation}
The fact is  that the compleatness  of the root vectors system   is not sufficient for the series convergence. In accordance with the definition given by Lidskii V.B. \cite{firstab_lit:1Lidskii}
 series \eqref{10y} is said to be summable  to the element $f$ via the method $(A,\lambda,s)$   if the following relation holds
\begin{equation}\label{11y}
\exists\{M_{\mu}\}_{0}^{\infty}\subset \mathbb{N}: \sum_{\mu=0}^{\infty}\sum_{q=M_{\mu}+1}^{ M_{\mu+1}}\mathcal{P}_{q}(s,t)f=S(t)f,
\end{equation}
$$
S(t)f\stackrel{\mathfrak{H}\,\;}{\rightarrow} f,\;t\rightarrow 0.
$$

\begin{deff}\label{D1} Assume that
 $$
 \exists\{M_{\mu}\}_{0}^{\infty}\subset \mathbb{N}: \sum_{\mu=0}^{\infty}\sum_{q=M_{\mu}+1}^{ M_{\mu+1}}\mathcal{P}_{q}(s,t)f=S(t)f,\,f\in \mathfrak{M}\subset \mathfrak{H},
$$
$$
S(t)f\stackrel{\mathfrak{H}\,\;}{\rightarrow} f,\;t\rightarrow 0,
$$
then the operator $B$ is said to be in the class $\mathcal{A}(s,\mathfrak{M}),$   i.e. in symbol    $B\in \mathcal{A}(s,\mathfrak{M}).$ The parameter $s$ is called by   the summation order.
\end{deff}

\section{ Main results}

\subsection{Splitting  to the infinite set of the invariant subspaces }

 Denote by $\mathfrak{M}_{k}$   the closure of the linear subspace of the root vectors  corresponding to an arbitrary subset of the eigenvalues  $\{\mu_{k_{j}}\}_{1}^{\infty}\subset \{\mu_{j}\}_{1}^{\infty}$ of the compact operator $B.$

\begin{lem}\label{L6}A compact operator    $B$  induces a  compact restriction $B_{k}$  on the invariant subspace $\mathfrak{M}_{k},$  moreover $$\Sigma(B_{k})=\{\mu_{k_{j}}\}_{1}^{\infty}.$$
\end{lem}
\begin{proof}

Let us show that the subspace $\mathfrak{M}_{k}$ is an invariant subspace of the operator $B.$     It is clear that the operator $B$ preserves linear combinations of the root vectors since $\mathfrak{N}_{q},\,q\in \mathbb{N}$ are invariant subspaces of the operator $B,$ see formula \eqref{6}. Thus, it suffices to show that the   images of  the  limits of the root vectors linear combinations belong to $ \mathfrak{M}_{k}.$   To prove the fact, consider an element $g$ such that
$$
Bf_{n}\stackrel{\mathfrak{H}\,\;}{\rightarrow} g ,\, f_{n}\stackrel{\mathfrak{H}\,\;}{\rightarrow} f \in \,\mathfrak{M}_{k},\,n \rightarrow \infty,\, f_{n}:=\sum\limits_{\nu=0}^{n}e_{\nu}c_{n\nu},
$$
where $c_{n\nu}$ are complex valued coefficients, $e_{\nu}$  root vectors corresponding to the set $\{\mu_{k_{j}}\}_{1}^{\infty}.$ In accordance with the above, we have $Bf_{n}\in \mathfrak{M}_{k}.$ Since  $ \mathfrak{M}_{k}$ is a closed subspace in the sense of the norm of the  Hilbert space $\mathfrak{H},$ then $g\in \mathfrak{M}_{k}.$ Therefore, $\mathfrak{M}_{k}$ is an invariant subspace of the operator $B.$ Note  that the restriction  $B_{k}$ is compact, since $B$ is compact.

  Let us prove the fact $\Sigma(B_{k})=\{\mu_{k_{j}}\}_{1}^{\infty}.$  Note that  in accordance with the Hilbert theorem the  spectrum of a compact operator  except for the point zero consists of normal  eigenvalues. Thus, it suffices to prove that the set of the  eigenvalues of the operator $B_{k}$ coincides  with the set $\{\mu_{k_{j}}\}_{1}^{\infty}.$
  Consider the set of the eigenvalues $\{\mu_{n_{j}}\}_{1}^{\infty}=\{\mu_{j}\}_{1}^{\infty}\setminus \{\mu_{k_{j}}\}_{1}^{\infty},$  then  in accordance with Theorem 6.17, Chapter III \cite{firstab_lit:kato1980},    we  have the decomposition
\begin{equation*}
\mathfrak{H}=  \mathfrak{M}_{l}'\oplus  \mathfrak{M}_{l}'',\,l\in \mathbb{N},
\end{equation*}
corresponding to the finite set $\{\mu_{n_{j}}\}_{1}^{l},$
where $ \mathfrak{M}_{l}'$ is   a finite dimensional invariant subspace of the operator $B$ generated by the root vectors corresponding to $\{\mu_{n_{j}}\}_{1}^{l},$  and  $\mathfrak{M}_{l}''$ is its parallel complement respectively, we have   $P_{l}\mathfrak{H}=\mathfrak{M}_{l}',\;(I-P_{l})\mathfrak{H}=\mathfrak{M}_{l}'',$ where $P_{l}$ is the corresponding projector, i.e.
$$
 P_{ l }f:= -\frac{1}{2\pi i}\oint_{\Gamma_{l}}(B-\lambda I)^{-1}f\,d\lambda,\; f\in \mathfrak{H},
$$
the contour  $\Gamma_{l}$ is a closed contour bounding a domain containing the set of the  eigenvalues  $\{\mu_{n_{j}}\}_{1}^{l}$ only.
Observe that the operator $P_{l}$ is bounded in the Hilbert space  $\mathfrak{H}.$ It can be proved easily since the space $\mathfrak{M}_{l}'$ is finite dimensional, thus  using  the orthogonalization procedure in the Hilbert space $\mathfrak{H}$ having preserved the basis vectors belonging to $ \mathfrak{M}_{l}'$  we easily obtain
$$
(P_{l}f,P_{l}f)_{\mathfrak{H}}\leq(f,f)_{\mathfrak{H}},
$$
from what follows that $\|P_{l}\|\leq 1.$ Let us show that the subspace $\mathfrak{M}_{l}''$ is closed, assume that  $$
 g_{k}\stackrel{\mathfrak{H}\,\;}{\rightarrow} g,\,k\rightarrow\infty,\,\{g_{k}\}_{1}^{\infty}\subset \mathfrak{M}_{l}'',
 $$
in accordance with the continuous property of the operator $P_{l},$ taking into account $P_{l}g_{k}=0,$ we get
$$
\|P_{l}g\|_{\mathfrak{H}}=\|P_{l}(g_{k}-g)\|_{\mathfrak{H}}\leq \| g_{k}-g \|_{\mathfrak{H}},\Rightarrow P_{l}g=0,\Rightarrow g\in \mathfrak{M}_{l}''.
$$
Thus, we conclude that the space $\mathfrak{M}_{l}''$ is closed. Note  that in accordance with Theorem 6.17, Chapter III \cite{firstab_lit:kato1980}, we have  $P_{l}e=0,$ where
 $e$ is a root vector corresponding to the eigenvalue $\mu\in \{\mu_{j}\}_{1}^{\infty}\setminus\{\mu_{n_{j}}\}_{1}^{l}.$ Hence    the closure of the root vectors linear combinations corresponding to $\{\mu_{j}\}_{1}^{\infty}\setminus\{\mu_{n_{j}}\}_{1}^{l}$  belongs to $\mathfrak{M}_{l}''.$   Therefore
   $\mathfrak{M}_{k}\subset\mathfrak{M}_{l}'',\,l\in \mathbb{N},$ since $\{\mu_{j}\}_{1}^{\infty}\setminus\{\mu_{n_{j}}\}_{1}^{l}\supset \{\mu_{k_{j}}\}_{1}^{\infty}.$

 Note that  in accordance with the   made assumptions   the root vectors system corresponding to the set of the eigenvalues $\{\mu_{k_{j}}\}_{1}^{\infty} $ belongs to
the root vectors system of the operator  $B_{k}.$    Let us show that they are coincided, i.e. there does not exist a root vector of the operator $B_{k}$ corresponding to an eigenvalue that differs from $\{\mu_{k_{j}}\}_{1}^{\infty} .$  Assume the contrary, then taking into account the fact $B_{k}\subset B,$   we should admit that there exists a number $N$ and an eigenvalue
$
\mu\in \{\mu_{n_{j}}\}_{1}^{N},
$
so that $(B-\mu I)^{\xi}e=0,\,e \in \mathfrak{M}_{k},\,\xi\in \mathbb{N}.$ Hence  $\mathfrak{M}_{k}\cap \mathfrak{M}_{p}' \neq0,\,p\geq N$ but it contradicts the proved above fact  in accordance with  which $\mathfrak{M}_{k}\subset\mathfrak{M}_{p}'',$ since $\mathfrak{M}'_{p}\cap\mathfrak{M}_{p}''=0.$  Therefore the root vectors system of the operator $B_{k}$ coincides with the root vectors system  of the operator $B$ corresponding to the set of the eigenvalues $\{\mu_{k_{j}}\}_{1}^{\infty} .$ It implies that  the set of the  eigenvalues of the operator $B_{k}$ coincides  with the set $\{\mu_{k_{j}}\}_{1}^{\infty}.$  The proof is complete.

\end{proof}

\subsection{Splitting of the counting function}\label{P3.2}

Consider a subsequence of the natural numbers
\begin{equation}\label{15}
 N_{\nu}= \sum\limits_{k=0}^{\nu}[\nu^{\beta}-k^{\beta}],\,\beta>0,\,\nu\in \mathbb{N}_{0}.
 \end{equation}
Let us   split   the sequence of the principal characteristic numbers    $\{z_{j}\}_{1}^{\infty}$  on the   groups $\{z_{k_{j}}\}_{1}^{\infty},$ i.e.
\begin{equation}\label{16}
\{z_{j}\}_{1}^{\infty}=\bigcup\limits_{k=0}^{\infty} \{z_{k_{j}}\}_{1}^{\infty},
\end{equation}
corresponding to the numbers $N_{k\nu}:=[\nu^{\beta}-k^{\beta}];\,N_{k\nu}=0, \nu\leqslant k$  so that   the disk  $\{z:\, |z|\leq|z_{N_{\nu}}|\},$ where we formally put $z_{0}:=0,$  contains $N_{k\nu}$ elements of the $k$-th group $\{z_{k_{j}}\}_{1}^{\infty}.$
In terms of counting functions, we have
$$
n(|z_{N_{\nu}}|,z_{ j })=N_{\nu}=\sum\limits_{k=0}^{\nu}[\nu^{\beta}-k^{\beta}]=\sum\limits_{k=0}^{\nu}n(|z_{N_{\nu}}|,z_{k_{j}}).
$$
Here, we ought to point out  that   $\{k_{j}\}_{1}^{\infty}$ is a subsequence of natural numbers defined by the index $k$ and in accordance with the last union, we have
 $$
\mathbb{N}=\bigcup _{k=0}^{\infty} \{k_{j}\}_{1}^{\infty}.
$$
It is clear that     splitting \eqref{16}  induces, in the natural way, the splitting
$$
\{\lambda_{j}\}_{1}^{\infty}=\bigcup\limits_{k=0}^{\infty} \{\lambda_{k_{j}}\}_{1}^{\infty}.
$$
In accordance with the above, we can express the principal  index
$$
p_{N_{\nu}}=\sum\limits_{k=0}^{\nu}\sum\limits_{j=1}^{[\nu^{\beta}-k^{\beta}]}\Delta_{j}(k),
$$
where $\Delta_{j}(k)$ denotes  algebraic multiplicity    corresponding to the  principal characteristic number
$z_{k_{j}}\in \{z_{k_{j}}\}_{1}^{\infty}.$
  For a convenient form of writing,   we will use the following shortages
$$
 \tilde{N}_{k\nu}:=\sum\limits_{j=1}^{[\nu^{\beta}-k^{\beta}]}\Delta_{j}(k),\;\tilde{N}_{\nu}:=p_{N_{\nu}}.
$$
It is rather  clear that
 the disk on  the  complex plane
$
\{z\in \mathbb{C}:\, |z|\leq|\lambda_{\tilde{N}_{\nu}}|\}
$
contains   $\tilde{N}_{k\nu}$ characteristic numbers belonging to    the  $k$-th group $\{\lambda_{k_{j}}\}_{1}^{\infty},$ where we put $\lambda_{0}:=0$  in correspondence with  the formalities accepted above.

 Further, applying Lemma \ref{L6}, we put the  operator  $B_{k}$ in correspondence with the $k$-th group $\{\lambda_{k_{j}}\}_{1}^{\infty}$ and  use the following  notation     $\lambda_{j}(B_{k}):=\lambda_{k_{j}}.$
\begin{deff}\label{D2}
Assume that
$$
\Delta_{j}<Cj^{\phi},\,0<\phi<1,\,j\in \mathbb{N},
$$
then we will say that   the operator $B$ has the sequence of the     algebraic multiplicities
    of the {$\phi$-th growth} (of the  lowest   growth if $\phi$ can be chosen arbitrary small).
\end{deff}
\begin{lem}\label{L7} Assume that $B\in \mathfrak{S}_{\sigma},\,0<\sigma<\infty,$ has  the sequence of the  multiplicities     of the  $\phi$-th growth,
then
\begin{equation*}
\lim\limits_{r\rightarrow\infty} \frac{n(r,B_{k})}{r^{  s}}=0, \,s> \sigma\left(\frac{\beta}{\beta+1}+\phi\right),
\end{equation*}
  uniformly  with respect to $k\in \mathbb{N}_{0}.$
\end{lem}
\begin{proof}

Consider  a subsequence of the natural numbers $\{ N_{\nu}\}^{\infty}_{0}$   defined in \eqref{15}. Let us   prove the following asymptotic  formula
\begin{equation}\label{17}
N_{\nu}\sim \frac{\beta}{\gamma}\,\nu^{\gamma},\;\nu\rightarrow\infty,
\end{equation}
here and further  $\gamma:=\beta+1.$ For this purpose, we will  estimate the given sum  by  a corresponding  definite integral, i.e. calculating the integral, we have on the one hand
$$
\sum\limits_{k=1}^{\nu}k^{ \beta}\leq\int\limits_{1}^{\nu+1}x^{\beta}dx=\frac{(\nu+1)^{\gamma}}{\gamma}-\frac{1}{\gamma },
$$
on the other hand
$$
\sum\limits_{k=1}^{\nu}k^{ \beta }=\sum\limits_{k=2}^{\nu}k^{ \beta }+1\geq \int\limits_{1}^{\nu }x^{\beta }dx+1=\frac{ \nu^{\gamma }}{\gamma}+\frac{\beta}{\gamma}.
$$
Therefore
$$
N_{\nu}=\sum\limits_{k=0}^{\nu-1}[\nu^{\beta}-k^{\beta}]\geq \nu^{\gamma} -\nu-\frac{ \nu ^{\gamma}}{\gamma}+\frac{1}{\gamma }=\frac{ \beta\nu ^{\gamma}}{\gamma} -\nu+ \frac{1}{\gamma }.
$$
Analogously
$$
N_{\nu}=\sum\limits_{k=0}^{\nu}[\nu^{\beta}-k^{\beta}]\leq \sum\limits_{k=0}^{\nu} (\nu^{\beta}-k^{\beta}) \leq \nu^{\beta}(\nu+1)-\frac{ \nu^{\gamma }}{\gamma}-\frac{\beta}{\gamma}= \frac{ \beta\nu ^{\gamma}}{\gamma}+\nu^{\beta}-\frac{\beta}{\gamma},
$$
from what follows the desired result. Note that in accordance with the fact that the operator belongs to the Schatten-von Neumann  class $\mathfrak{S}_{\sigma},$  we have
$$
\lim\limits_{r\rightarrow\infty} \frac{n(r,B)}{r^{\sigma} }=0.
$$
This fact obviously follows   from the implication
\begin{equation}\label{18}
 B\in \mathfrak{S}_{\sigma}   \Rightarrow s_{n}(B)=o(n^{-1/\sigma}) \Rightarrow |\mu_{n}(B)|=o(n^{-1/\sigma}),\,n\rightarrow\infty,\, 0<\sigma<\infty,
\end{equation}
  see $8^{\circ},$ \S7,  Chapter III, \cite{firstab_lit:1Gohberg1965},  Corollary 3.2,  \S3,  Chapter II,   \cite{firstab_lit:1Gohberg1965}. Note  that
  $
  \Delta_{j}\geq1,
  $
hence
$$
 \tilde{N}_{\nu} =\sum\limits_{k=0}^{\nu}\sum\limits_{j=1}^{[\nu^{\beta}-k^{\beta}]}\Delta_{j}(k)\geq \sum\limits_{k=0}^{\nu}[\nu^{\beta}-k^{\beta}]=N_{\nu}.
$$
Therefore, applying  asymptotic  formula \eqref{17},   we obtain
\begin{equation}\label{19}
\frac{\nu^{\gamma}  }{|\lambda _{\tilde{N}_{\nu} }|^{\sigma}}\leq C  \frac{\tilde{N} _{\nu } }{|\lambda _{\tilde{N}_{\nu} }|^{\sigma}} \leq C_{\nu},\;C_{\nu}\rightarrow 0,\,\nu\rightarrow\infty.
 \end{equation}
  Note that in accordance with the splitting, assuming that $N_{k\nu}>0,$ we have
 $$
 z_{k_{j}}\in \{z\in \mathbb{C}:\,|z|\leq z_{N_{\nu}}\},\,j=1,2,...,N_{k\nu},
 $$
 therefore using the condition imposed upon the growth of the algebraic multiplicities, we obtain
 $$
 \Delta_{j}(k)\leq \max\limits_{j\in [1, \,N_{\nu}]}\Delta_{j}\leq C N_{\nu}^{\phi}.
 $$ 
Applying   formula \eqref{17}, we get
\begin{equation}\label{20}
\tilde{N}_{k\nu  }=\sum\limits_{j=1}^{[\nu^{\beta}-k^{\beta}]}\Delta_{j}(k)\leq \sum\limits_{j=1}^{[\nu^{\beta}-k^{\beta}]}C N^{\phi}_{\nu}=
 C  [\nu^{\beta}-k^{\beta}] N^{\phi}_{\nu}\leq C \nu^{\,\beta+\phi\gamma}.
\end{equation}
Assuming that $k$ is fixed, consider a sequence
$$
 \frac{\tilde{N}_{k\nu  }}{ |\lambda _{\tilde{N}_{k\nu } }(B_{k}) |^{s }},\,s>\sigma\left(\frac{\beta}{\gamma}+\phi\right), \,\nu =k+1,k+2,...\,.
$$
Observe that the numbers $\tilde{N}_{k\nu  }$    have variations only corresponding to the the values of the index $\nu$ satisfying the condition $\left[(\nu+1)^{\beta}-k^{\beta}\right]>[\nu ^{\beta}-k^{\beta}],$ in this case, we have
 $$
|\lambda _{\tilde{N}_{ \nu  }}|<|\lambda _{\tilde{N}_{k\nu}+q}(B_{k})|\leq |\lambda _{\tilde{N}_{\nu +1}}|,\,0 < q\leq\tilde{N}_{k\nu+1}- \tilde{N}_{k\nu }.
$$
Using the lower estimate, applying  \eqref{20},  we get
$$
\frac{\tilde{N}_{k\nu}+q}{ |\lambda _{\tilde{N}_{k\nu}+q }(B_{k}) |^{s  }}<
\frac{\tilde{N}_{k\nu+1  }}{  |\lambda  _{\tilde{N}_{\nu } }  |^{s  }}\leq C
 \left\{\frac{\nu^{\gamma}  }{|\lambda _{\tilde{N}_{\nu} }|^{\sigma}}\right\}^{\frac{\beta}{\gamma}+\phi}\leq C C^{\frac{\beta}{\gamma}+\phi}_{\nu},
$$
from what follows the desired result. The proof is complete.
\end{proof}

  Using the latter  result, we can represent   a scheme of reasonings allowing to    decrease the summation  order. The following   paragraph is devoted to a sharper estimate for the canonical product, however the final aim is to improve  the estimate for the norm of the resolvent what can be implemented
due to  the properties  of the given above  artificially constructed subsequence of the eigenvalues. Since the result is fundamental and relates  to the issue in the framework of the infinite determinant theory   Chapter IV  \cite{firstab_lit:1Gohberg1965}, we may claim that  it may represent the interest itself.    \\

\subsection{  Sharper estimate for the canonical product}

 \begin{lem}\label{L8} Assume that $B$ is  a compact  operator $B\in \mathfrak{S}_{1},\,\Theta(B)\subset \mathfrak{L}_{0}(\theta),\,\theta<\pi/4,$     then the following estimate holds
$$
\prod\limits_{n=1}^{\infty}|1-\lambda\mu_{n}(Q_{1}BQ_{1})| \leq \prod\limits_{n=1}^{\infty}|1+\lambda\mu_{n}( B )|,\,|\mathrm{arg}\lambda|<\pi/4,
$$
where $Q_{1}$ is the orthogonal projector corresponding to the orthogonal complement of the one-dimensional  subspace generated by   an element $f\in \mathfrak{H}.$
\end{lem}
\begin{proof}
Firstly, we should note that
in accordance with   Lemma 1 \cite{firstab_lit:1Lidskii}, we have
\begin{equation*}
s_{n}(Q_{1}BQ_{1} )\leq s_{n}( B ),\,n=1,2,...\,,
\end{equation*}
hence $Q_{1}BQ_{1}\in \mathfrak{S}_{1}.$  Note that by virtue of the relation
$$
   \mathrm{Re} (Q_{1}BQ_{1}f,f)=\mathrm{Re}(BQ_{1}f,Q_{1}f)  \geq 0 ,\,f\in \mathfrak{H},
$$
we obtain the fact $\Theta(Q_{1}BQ_{1})\subset \Theta( B).$
  Consider the operators $ B(\lambda)=(I+\lambda  B )$ and $ B_{1}(\lambda):=(Q_{1}+\lambda  B_{1} ),\,B_{1}:=Q_{1}BQ_{1}.$
Note that
 $$
 B^{\ast}(\lambda)B(\lambda)=I+C(\lambda),\,C(\lambda)= |\lambda B|^{2}+2\mathfrak{Re} (\lambda B),\;
 $$
 $$
 B^{\ast}_{1}(\lambda)B_{1}(\lambda)=Q_{1}+C_{1}(\lambda),\,C_{1}(\lambda):= |\lambda B_{1}|^{2}+2\mathfrak{Re} (\lambda B_{1}),
 $$
 where $|B|^{2}:=B^{\ast}B.$ It is clear that  since $C(\lambda)$ is compact selfadjoint then the set of the eigenvectors   is complete in $\overline{\mathrm{R}(C(\lambda))},$  and  we can choose a basis $\{e_{k}\}_{1}^{\infty}$ in $\overline{\mathrm{R}(C(\lambda))}$   such that the operator  matrix will have a diagonal form - the eigenvalues are situated on the major diagonal. It is clear that the same reasonings are true for the operator $C_{1}(\lambda).$ Applying Corollary 2.2, \S2, Chapter II, \cite{firstab_lit:1Gohberg1965}, we obtain easily  the fact    $C(\lambda),C_{1}(\lambda)\in \mathfrak{S}_{1},$
  therefore,   applying Corollary 1.1, $ 3^{\circ},$  Chapter IV, \S1, \cite{firstab_lit:1Gohberg1965}, we get
 $$
  \lim\limits_{n\rightarrow\infty}\det\{P_{n}+P_{n}C(\lambda)P_{n}\}=\det\{I+ C(\lambda) \},
 $$
 where $P_{n},\,n\in \mathbb{N}$ is an orthogonal projector into the subspace generated by the eigenvectors $\{e_{k}\}_{1}^{n}$   of the operator $C(\lambda).$ Note that
$
 Q_{1}P_{1n} = P_{1n} ,
$
where $P_{1n}$ is an orthogonal projector into the subspace generated by the eigenvectors $\{e_{1k}\}_{1}^{n}$ of the operator $C_{1} (\lambda).$
   Analogously, we get
$$
 \lim\limits_{n\rightarrow\infty}\det\{P_{1n}+P_{1n}C_{1}(\lambda)P_{1n}\}=\det\{I+ C_{1}(\lambda) \}.
$$
Consider
$$
(C_{1}(\lambda)f,f)=|\lambda|^{2}(Q_{1}B^{\ast}Q_{1} B Q_{1}f,f)+2(\mathfrak{Re} (\lambda B_{1})f,f)\leq |\lambda|^{2}(B^{\ast}  B Q_{1}f,Q_{1}f)+2(\mathfrak{Re} (\lambda B )Q_{1}f,Q_{1}f)=
$$
$$
=(Q_{1}C (\lambda)Q_{1}f,f),
$$
here we used the obvious relations
$$
(Q_{1}B^{\ast}Q_{1} B Q_{1}f,f)\leq (B^{\ast} B Q_{1}f,Q_{1}f),\;\mathfrak{Re} (\lambda B_{1})=Q_{1}\mathfrak{Re} (\lambda B)Q_{1}.
$$
Since  $C (\lambda),\, |\arg  \lambda|<\pi/4$ is a compact non-negative selfadjoint operator, then by virtue  to the minimax principle for the eigenvalues  (see Courant theorem  \cite[p.120]{firstab_lit: Courant-Hilbert}), we get
$$
\mu_{n}(C_{1} (\lambda))\leq\mu_{n}( Q_{1}C (\lambda)Q_{1})\leq \mu_{n}(C (\lambda)),\,n\in \mathbb{N},
$$
therefore
$$
\det\{P_{1n}+P_{1n}C_{1}(\lambda)P_{1n}\}\leq\det\{P_{n}+P_{n}C(\lambda)P_{n}\},\,n\in \mathbb{N}.
$$
Passing to the limit, we get
\begin{equation}\label{21}
\det\{I+C_{1}(\lambda) \}\leq\det\{I+C(\lambda) \}.
\end{equation}
On the other hand,  in accordance with  Theorem 2.3 Chapter V  \cite{firstab_lit:1Gohberg1965} the system of the root vectors (including the root vectors corresponding to the zero eigenvalue) of the operator $iB$ is compleat in $\mathfrak{H}.$   Indeed, we will prove it if we show that $iB$ is dissipative and  $iB\in \mathfrak{S}_{1}.$  Taking into account the fact  $\Theta(B) \subset \mathfrak{L}_{0}(\theta),\,\theta<\pi/4,$   we conclude that
$$
\mathrm{Im}(iBf,f)=\mathrm{Re}(Bf,f)\geq 0,\,f\in \mathfrak{H}.
$$
Therefore, the operator $iB$ is dissipative. It is clear that $s_{n} (iB)=s_{n} (B),\,n\in \mathbb{N},$ hence $iB\in \mathfrak{S}_{1}.$ Thus, we obtain the desired result.
Note that in accordance with   Theorem 2.3   Chapter V  \cite{firstab_lit:1Gohberg1965}, we have
$$
\mathfrak{H}=\mathfrak{C}_{0}(iB)\, \dot{+ }\,  \mathfrak{H}_{0} (iB),
$$
where $\mathfrak{C}_{0}(iB)$ is the invariant subspace generated by the closure of the linear combinations of the root vectors corresponding to   non-real  eigenvalues of the operator $iB$ and $\mathfrak{H}_{0} (iB)$ is the invariant subspace on which the restriction of  $iB$ is selfadjoint. Since $\Theta(B)\subset \mathfrak{L}_{0}(\theta)$ then $\mathfrak{C}_{0}(iB)=\mathfrak{C}(iB)=\mathfrak{C}(B)$ the latter symbol denotes the closure of the linear combinations of the root vectors of the operator $B,$ we used the fact that the operators $B$ and $iB$ have the same root vectors.  It implies that  $\mathfrak{H}_{0} (iB)=\mathrm{N}(B),$ since the operator $iB$ does not have real eigenvalues.  Hence
$$
\mathfrak{H}=\mathfrak{C}(B)\, \dot{+ }\, \mathrm{N}(B).
$$
Therefore, in accordance with Lemma 4.1 Chapter I \cite{firstab_lit:1Gohberg1965}, we can construct   an orthogonal Schur basis $\{\omega_{j}\}_{1}^{\infty}\subset \mathfrak{C} $   so that the  matrix
of the operator induced in  $\mathfrak{C}$   has a triangle form. Choosing an arbitrary basis in the space   $\mathrm{N}(B),$
 uniting   bases of the orthogonal decomposition, we obtain the fact that the matrix of the  operator $B$  has a triangle form in a newly constructed united basis.
 Thus, choosing   orthogonal projectors   corresponding to n-dimensional subspaces, the property of the triangle determinant, we have
$$
\det\{P_{n} B (\lambda)P_{n}\}=\prod\limits_{k=1}^{n}\{1+\lambda\mu_{k}( B )\};\;\det\{P_{n} B^{\ast} (\lambda)P_{n}\}=\det\{[P_{n} B (\lambda)P_{n}]^{\ast}\}=\prod\limits_{k=1}^{n}\{1+\overline{\lambda\mu_{k}( B )}\},
$$
therefore
$$
\det\{P_{n}+P_{n}C(\lambda)P_{n}\}=\det\{P_{n} B (\lambda)P_{n}\}\det\{P_{n} B^{\ast} (\lambda)P_{n}\}=   \prod\limits_{k=1}^{n}|1+\lambda\mu_{k}( B )|^{2}.
$$
Analogously, we get
$$
\det\{P_{1n}+P_{1n}C_{1}(\lambda)P_{1n}\}=\det\{P_{1n} B_{1} (\lambda)P_{1n}\}\det\{P_{1n} B^{\ast}_{1} (\lambda)P_{1n}\}=\prod\limits_{k=1}^{n}|1+\lambda\mu_{k}(Q_{1}BQ_{1})|^{2}.
$$
Applying Corollary 1.1, $ 2^{\circ},$  Chapter IV, \S1, \cite{firstab_lit:1Gohberg1965}, we conclude that
$$
 \lim\limits_{n\rightarrow\infty}\det\{P_{n}+P_{n}C(\lambda)P_{n}\} =\det\{I+C(\lambda) \},\;\lim\limits_{n\rightarrow\infty}\det\{P_{1n}+P_{1n}C_{1}(\lambda)P_{1n}\} =\det\{I+ C_{1}(\lambda) \}.
$$
Taking into account\eqref{21}, we get
$$
\prod\limits_{n=1}^{\infty}|1+\lambda\mu_{n}(Q_{1}BQ_{1})| \leq \prod\limits_{n=1}^{\infty}|1+\lambda\mu_{n}( B )|.
$$
Having noticed the fact
$$
|1-\lambda\mu_{n}(Q_{1}BQ_{1})| \leq  |1+\lambda\mu_{n}(Q_{1}BQ_{1})| ,\,|\arg \lambda|<\pi/4,
$$
we obtain the desired result. The proof is complete.
\end{proof}

\begin{lem}\label{L9}   Assume that    $B\in \mathfrak{S}_{ \sigma},\,0<\sigma\leq1, \,\Theta(B)\in \mathfrak{L}_{0}(\theta),\,\theta<\pi/4,$  the following relation holds
$$
\lim \limits_{r\rightarrow\infty}\frac{n(r,B)}{r^{s}}=0,\;0<c<s<\sigma.
$$
Then
for arbitrary numbers  $R,\delta$   such that $R>0,\,0<\delta<1,$ there exists  $(1-\delta)R<\tilde{R}<R,$
   so that the following estimate holds
$$
 \|(I-\lambda B )^{-1}\|\leq   2e^{h\left(  2e R \right)\ln\frac{45e^{4}}{\delta}}   ,\,|\lambda|=\tilde{R},\,|\arg \lambda|<\pi/4,
$$
 where
$$
 h(r )=   \left(\int\limits_{0}^{r}\frac{n (t,B)dt}{t }+
r \int\limits_{r}^{\infty}\frac{n (t,B)dt}{t^{ 2  }}\right).
$$
\end{lem}
\begin{proof}
 Consider a Fredholm  determinant of the operator $B,$ in accordance with Lemma \ref{L3} it has a representation
$$
\Delta(\lambda)= \prod\limits_{n=1}^{\infty}\left\{1- \lambda \mu_{n}(B)    \right\},\,\lambda\in \mathbb{C}.
$$
 Let us chose an arbitrary element $f\in \mathfrak{H}$  and construct a new orthonormal basis having put $f$ as a first basis element. Note that   relations \eqref{4} hold for the matrix coefficients of the operator $B$ in a new basis, this fact follows from the well-known theorem for the operator class $\mathfrak{S}_{1}.$  Thus, using the given  above representation for the resolvent \eqref{5}, we obtain
\begin{equation}\label{18b}
\Delta (\lambda)\left((I-\lambda B )^{-1}f,f\right)_{\mathfrak{H}}= \Delta^{11}(\lambda),
\end{equation}
the latter entire function  depends on the choice  of an  element $f\in\mathfrak{ H}$ and   we   reflect   this fact in the  notation $\tilde{f}(\lambda):=\Delta^{11}(\lambda).$ Let us notice  the fact that  $\tilde{f}(\lambda)$  represents   the Fredholm determinant of the operator $Q_{1}B Q_{1},$ where $Q_{1}$  is the projector into  orthogonal  complement of the element $f\in \mathfrak{H}.$    Having applied  Lemma \ref{L2}  (Lemma 1 \cite{firstab_lit:1Lidskii}), we obtain
$$
s_{n}(Q_{1}B Q_{1})\leq s_{n}( B  ),\,n\in \mathbb{N}.
$$
Applying Lemma \ref{L3}, we obtain  the representation
$$
\tilde{f}(\lambda)=\prod\limits_{n=1}^{\infty}\left\{1- \lambda \mu_{n}(Q_{1}B Q_{1})    \right\},\,\lambda\in \mathbb{C}.
$$
Applying Lemma \ref{L8}, we have
\begin{equation}\label{22a}
|\tilde{f}(\lambda)|=\prod\limits_{n=1}^{\infty}\left|1-  \lambda \mu_{n}(Q_{1}B Q_{1})\right|\leq \prod\limits_{n=1}^{\infty}|1+\lambda\mu_{n}( B  )|  \leq \prod\limits_{n=1}^{\infty}\left\{1+ |\lambda \mu_{n}( B  )| \right\},\,
$$
$$
f\in \mathfrak{H},\,|\arg \lambda|<\pi/4.
\end{equation}
 Let us  prove   the following relation
\begin{equation}\label{22}
|\Delta (\lambda)|\cdot\|(I-\lambda B )^{-1}\|\leq 2 \prod\limits_{n=1}^{\infty}\left\{1+ |\lambda \mu_{n}( B  )| \right\},\,|\arg \lambda|<\pi/4.
\end{equation}
For this purpose, define an  operator
$
D_{B}(\lambda):=\Delta (\lambda) (I-\lambda B )^{-1}
$
 in accordance with \eqref{18b}, we have a correspondence between the notations   $(D_{B } (\lambda)f,f )_{\mathfrak{H}}=\tilde{f}(\lambda).$
 Consider the   decomposition on the    Hermitian components
$$
D_{B } (\lambda)=\mathfrak{Re}D_{B } (\lambda)+i \,\mathfrak{Im}D_{B } (\lambda).
$$
Note that the  Hermitian components are selfadjoint operators.
Thus, using the well-known  formula for the norm of a selfadjoint operator,  for an arbitrary fixed $\lambda\in \mathbb{C},$ we get
$$
\|D_{B }(\lambda)\|=\sup \limits_{\|f\|  \leq 1 }\|\mathfrak{Re} D_{B }(\lambda)f+i\, \mathfrak{Im} D_{B }(\lambda)f \|_{\mathfrak{H}} \leq\sup \limits_{\|f\|\leq 1 }\|\mathfrak{Re} D_{B }(\lambda)f  \|_{\mathfrak{H}}+ \sup \limits_{\|f\|\leq 1 }\|  \mathfrak{Im} D_{B }(\lambda)f \|_{\mathfrak{H}}=
$$
$$
=\sup \limits_{\|f\| =1 }|(\mathfrak{Re} D_{B } (\lambda)f,f )_{\mathfrak{H}} |+ \sup \limits_{\|f\| =1 }| (\mathfrak{Im} D_{B } (\lambda)f,f )_{\mathfrak{H}}|=
 $$
 $$
 =\sup \limits_{\|f\| =1 }|\mathrm{Re}(  D_{B } (\lambda)f,f )_{\mathfrak{H}}  |+ \sup \limits_{\|f\| =1 }|\mathrm{Im}(  D_{B } (\lambda)f,f )_{\mathfrak{H}}  |\leq 2\sup \limits_{\|f\| =1 }| (  D_{B } (\lambda)f,f )_{\mathfrak{H}}  |= 2\sup \limits_{\|f\|= 1 } |  \tilde{f}(\lambda) |.
$$
Taking into account \eqref{22a}, we obtain \eqref{22}. 
In accordance with the made  assumptions,   we have   $n(r,B)=o( r^{s}),\,0<c<s<1,$ hence
$$
\sum\limits_{n=1}^{\infty} |\mu_{n}( B  )|  <\infty,
$$
therefore, applying  Lemma \ref{L1}  to the canonical product, we obtain
$$
|\Delta (\lambda)|\cdot\|(I-\lambda B )^{-1}\|\leq 2\prod\limits_{n=1}^{\infty}\left\{1+ |\lambda \mu_{n}( B  )| \right\}\leq  2e^{h (|\lambda|) },\,|\arg \lambda|<\pi/4.
$$
Now, to obtain   the lemma statement it suffices to estimate the absolute value of the  Fredholm determinant of the operator $B$ from below. For this purpose, let us notice that in accordance with Lemma \ref{L3} it is an entire function represented by the formula
$$
\Delta (\lambda)= \prod\limits_{n=1}^{\infty}\left\{1- \lambda \mu_{n}(B )    \right\},\,\lambda\in \mathbb{C}.
$$
In accordance with the  Joseph Cartan concept, we can obtain an  estimate from   below for the entire function that  holds on the complex plane except may be an   exceptional set of circulus. The latter cannot be found but and we are compelled to make an exclusion evaluating the measures. However, in the paper \cite{firstab_lit(frac2023)},  we   produce the method allowing to find exceptional set of   circulus.  In the simplified case,   we can use  Theorem 4 \cite[p.79]{firstab_lit:Eb. LevinE} giving the lower bound of the absolute value for an analytic function in the disk, we have
$$
\ln|\Delta(\lambda)|\geq  -\ln \left\{\max\limits_{\psi\in[0,2\pi]}|\Delta ( 2e R e^{i\psi} )|\right\} \ln \frac{15e^{3}}{\eta},\,|\lambda|\leq R,
$$
except for the exceptional set of circles with the sum of radii less that $\eta R,$  where $\eta$ is an arbitrary small positive number.  Thus,  to find the desired circle $\lambda= e^{i\psi}\tilde{R}$ belonging to the ring, i.e. $R(1-\delta )<\tilde{R}<R,$   we have to   choose $\eta$ satisfying the inequality
$$
2\eta R<R-R(1-\delta )=\delta R;\;\eta< \delta/2.
$$
Hence, having chosen $\eta =\delta/3,$ we get
$$
\ln|\Delta(\lambda)|\geq  -\ln \left\{\max\limits_{\psi\in[0,2\pi]}|\Delta ( 2e R e^{i\psi} )|\right\}\ln \frac{45e^{3}}{\delta} ,\,|\lambda|= \tilde{R},
$$
Analogously to the above,  applying  Lemma \ref{L1} to the canonical product, we get
$$
|\Delta(\lambda)|= \left|\prod\limits_{n=1}^{\infty}\left\{1- \lambda \mu_{n}(B)    \right\}\right|\leq  e^{h (|\lambda|) }.
$$
Therefore
$$
\ln \left\{\max\limits_{\psi\in[0,2\pi]}|\Delta ( 2e R e^{i\psi} )|\right\}\leq h (2eR)  .
$$
Substituting, we get
$$
\ln|\Delta(\lambda)|\geq  - h \left( 2e R   \right)\ln  \frac{45e^{3}}{\delta}  ;\;
|\Delta(\lambda)|\geq  e^{- h \left( 2e R  \right)\ln\frac{45e^{3}}{\delta}  } ,\,|\lambda|= \tilde{R}.
$$
Note that
\begin{equation}\label{23}
 \frac{d h}{dr}=  \int\limits_{r}^{\infty}\frac{n (t,B)dt}{t^{ 2  }}>0.
\end{equation}
Hence, using the monotonous property of the function $h(r),$  we conclude that
$
h ( \tilde{R})  <h \left( 2e R   \right).
$
Combining the upper and the lower estimates, we obtain
$$
 \|(I-\lambda B )^{-1}\|=\frac{\|D_{B }(\lambda)\|}{|\Delta (\lambda)|} \leq  2e^{\left(1+\ln\frac{45e^{3}}{\delta}\right) h \left( 2e R  \right)  },\;|\lambda|= \tilde{R},\,|\arg \lambda|<\pi/4.
$$
The proof is complete.
\end{proof}

\newpage

\subsection{Infinitesimalness  of the summation order}

The following theorem is formulated in terms of   paragraph \ref{P3.2}.

\begin{teo} \label{T1}

 Assume that $B\in \mathfrak{S}_{\sigma},\,0<\sigma\leq1,\,\Theta(B) \subset   \mathfrak{L}_{0}(\theta) ,\,\theta< \pi/4,$  has   the sequence of the algebraic   multiplicities   of the    $\phi$ -   growth,
then
$$
B\in \mathcal{A}(s,\mathrm{R}(B)),\ \sigma \phi<s<\pi/2\theta.
$$
\end{teo}
\begin{proof}

Firstly, let us note that since $iB$ is a dissipative operator, then  in accordance with Theorem 2.3 paragraph 2  Chapter V \cite{firstab_lit:1Gohberg1965}, the system of the root vectors of the operator (including the root vectors corresponding to the zero eigenvelue) is complete in $\mathfrak{H}.$     Using   notations of paragraph  \ref{P3.2}, having chosen an arbitrary small  $\beta>0$ let us rearrange the sequence of the  characteristic numbers $\{\lambda_{n} \}_{1}^{\infty}$ of the operator $B$    in the   groups $\{\lambda_{k_{j}}\}_{1}^{\infty}$ in accordance with \eqref{16}.  Applying  Lemma \ref{L6}, we can put the  operators $B_{k}$  into   correspondence  with  the sequences  $\{\lambda_{k_{j}}\}_{1}^{\infty}.$  In accordance with Lemma \ref{L7},     we have
\begin{equation*}
\lim\limits_{r\rightarrow\infty} \frac{n(r,B_{k})}{r^{s}}=0 ,\;k\in \mathbb{N}_{0},\;s> \sigma\left(\frac{\beta}{\gamma}+\phi \right),\,\gamma:=\beta+1.
\end{equation*}
  Consider orhtogonal  projectors $P_{k}$ corresponding to  the invariant  subspaces  $\mathfrak{M}_{k}$ defined in Lemma  \ref{L6}.   It is clear that    $B_{k}=P_{k}BP_{k}.$  Applying Lemma  \ref{L2}, we get
$$
s_{n}(B_{k})\leq s_{n}(B),\,n\in \mathbb{N}.
$$
Therefore $B_{k}\in \mathfrak{S}_{1},\,k\in \mathbb{N}_{0}.$
Define a   contour in the complex plane
\begin{equation}\label{22b}
 \vartheta:= \left\{z\in \mathbb{C}:\;|z|=r_{0},\,|\mathrm{arg} z|\leq \theta_{0}\right\}\cup\left\{z\in \mathbb{C}:\;|z|>r_{0},\; |\mathrm{arg} z|=\theta_{0}\right\},
 $$
 $$
\,r_{0}:=|\lambda_{1}|-\varsigma,\,\theta_{0}:=\theta+\varsigma,
\end{equation}
where    $\varsigma$ is an arbitrary small positive fixed number.
 Now consider a sequence of the radii  $R_{\mu}=a\mu  +b,\,\mu\in \mathbb{N}_{0},\,a>b=r_{0}$ and define $\delta_{\mu}$ from the condition $R_{\mu}(1-\delta_{\mu})=a\mu,$ we have
$$
\delta_{\mu}=\frac{1}{1+\mu  a/b },\; R_{\mu+1}(1-\delta_{\mu+1}) >R_{\mu},
$$
additionally without loss of generality, we assume that  $\{z_{j}\}_{1}^{\infty}\cap \{R_{\mu}\}_{0}^{\infty}=\emptyset .$
   Applying Lemma \ref{L9}, we obtain the fact that in each case    there exists a sequence of contours $\{\tilde{R}_{k\mu}\}^{\infty}_{0},\; R_{\mu}(1-\delta_{\mu})<\tilde{R}_{k\mu}< R_{\mu}$ such that the   estimates  hold
$$
\|(I-\lambda B_{k} )^{-1}\|\leq 2 e^{  h_{k}(  2e R_{\mu})\ln\frac{45e^{4}}{\delta_{\mu}}  },\,|\lambda|=\tilde{R}_{k\mu},\,|\arg \lambda|<\pi/4, \,k\in \mathbb{N}_{0},
$$
where
$$
 h_{k}(r ):=    \int\limits_{0}^{r}\frac{n(t,B_{k})dt}{t }+
r \int\limits_{r}^{\infty}\frac{n(t,B_{k})dt}{t^{ 2  }} .
$$
 Additionally, we can assume that  $\{z_{j}\}_{1}^{\infty}\cap  \tilde{R}_{k\mu} =\emptyset,\,k,\mu \in \mathbb{N}_{0}$ by virtue of   arbitrariness in the  choice of $\tilde{R}_{k\mu}$ dictated by  Theorem 4 \cite[p.79]{firstab_lit:Eb. LevinE}.
Estimating, we get
\begin{equation}\label{24}
\|(I-\lambda B_{k} )^{-1}\|\leq 2 e^{  h_{k}(  2e R_{\mu})\ln [45e^{4}(1+\mu a/b)]  }\leq e^{ C_{0} h_{k}(  2e R_{\mu})\ln \mu },\,|\lambda|=\tilde{R}_{k\mu},\,|\arg \lambda|<\pi/4,\;k\in \mathbb{N}_{0}.
\end{equation}
  Let us prove that
\begin{equation}\label{25}
 \lim\limits_{r\rightarrow\infty}\frac{h_{k} (r)}{r^{s}}=0,\,r\rightarrow\infty,
\end{equation}
uniformly with respect to $k\in \mathbb{N}_{0}.$
Without loss of generality, we assume that $s<1,$  applying Lemma \ref{L7}, estimating the counting function  under the integrals, we get
$$
\forall \varepsilon>0,\,\exists M(\varepsilon): h_{k} (r)< \varepsilon\left\{\int\limits_{0}^{r}t^{s-1}dt+
r \int\limits_{r}^{\infty}t^{s-2}dt\right\}= \frac{\varepsilon r^{s} }{s(1-s)},\,r>M(\varepsilon),\,k\in \mathbb{N}_{0},
$$
what means the desired result \eqref{25}.  Define a subsequence  of the natural numbers $\{M_{\mu}\}_{0}^{\infty}$  as follows   a number $M_{\mu}$ indicates a quantity of principal  characteristic numbers of the operator $B$ belonging to the   open disk  with the radius $ R_{\mu},$ i.e.
$$
 M_{\mu}:=\mathrm{card}\{j\in \mathbb{N} :\,|z_{j}|<R_{\mu}\}.
$$
Then
the group
$
z_{M_{\mu}+1},\,z_{M_{\mu}+2},...,z_{M_{\mu+1}}
$  of the  principal  characteristic numbers   is inclosed in the closed  contour formed by the intersection of the contour   $\vartheta$ with the circulus having radii $ R_{\mu},R_{\mu+1}.$   Following to  the Lidskii V.B. \cite{firstab_lit:1Lidskii} results  consider the following   relation
$$
\frac{1}{2\pi i}\oint\limits_{\vartheta(R_{m+1}) }  e^{- \lambda^{s}  t}B \left(I-\lambda B \right)^{-1}f d\lambda=-\sum\limits_{\mu=0}^{m}
\sum\limits_{q=M_{\mu}+1}^{M_{\mu+1}} \mathcal{P}_{q}(s,t)f,\;f\in   \mathfrak{H},\;m\in \mathbb{N},
$$
where $\sigma\phi<s<\pi/2\theta,$
$$
\vartheta(R_{m+1}):=  \left\{z\in \mathbb{C}: |z|=r_{0}, |z|=R_{m+1}, |\mathrm{arg} z|\leq \theta_{0}\right\}\cup\left\{z\in \mathbb{C}: r_{0}<|z|<R_{m+1},  |\mathrm{arg} z|=\theta_{0}\right\}.
$$
Observe that the inner sum contains $M_{\mu+1}-M_{\mu}$ terms,   we have
$$
 M_{m+1}=\sum\limits_{\mu=0}^{m}M_{\mu+1}-M_{\mu}.
$$
The Lidskii V.B. idea is to pass to the limit in the last integral when $m$ tends to infinity and in this way to prove the series convergence. However, there are some obstacles in evaluating the norm of the integral, we should take a value of the summation order more then the index of the Schatten-von Neumann  class. At the same time the  heuristic reasonings lead us to the hypotheses that the   decay  of the    exponential function  has a surplus, therefore the latter can be replaced in the construction or at least the infinitesimal value of the order can be considered.
Consider a sum
$$
\sum\limits_{k=0}^{\psi(\mu)}\sum\limits_{j=M_{k\mu}+1}^{M_{k\mu+1}} \mathcal{P}_{k_{j}}(s,t)f,
$$
where the operators in the inner sum  correspond to the principal characteristic numbers $ z_{k_{j}},\,j=M_{k\mu}+1,M_{k\mu}+2,..., M_{k\mu+1},$ (see paragraph \ref{P3.2})
the  number $M_{k\mu}$ indicates a quantity of principal characteristic numbers of the operator $B_{k}$ belonging to the  open disk with the radius $\tilde{R}_{k\mu},$ i.e.
$$
M_{k\mu}:=\mathrm{card}\{j\in \mathrm{N },\,|z_{k_{j}}|<\tilde{R}_{k\mu}\}.
$$
The symbol
$
\psi(\mu):= \mathrm{card}\{k\in \mathbb{N}_{0}: |z_{k_{j}}|< \tilde{R}_{k\mu+1}\}
$
 denotes a function indicating a quantity of operators $B_{k}$ having characteristic numbers inside the circles  with the radii $\tilde{R}_{k\mu+1}.$ Observe the   representation
$$
\sum\limits_{k=0}^{\psi(m)}M_{km+1}= \sum\limits_{k=0}^{\psi(m)}\sum\limits_{\mu=0}^{m}\{M_{k\mu+1}-M_{k\mu}\}=
\sum\limits_{\mu=0}^{m}\sum\limits_{k=0}^{\psi(\mu)}\{M_{k\mu+1}-M_{k\mu}\},
$$
here  we have taken into account the fact $M_{k\mu+1}=0,\,k>\psi(\mu).$
Thus, if we analyze a mutual arrangement of the radii $R_{\mu},\tilde{R}_{k\mu},$  we   come to the fact that there exists a natural number   $m$ such  that
$$
 M_{m+1} \neq\sum\limits_{k=0}^{\psi(m)}M_{km+1}.
$$
The latter relation can be rewritten in terms of operators
$$
 \sum\limits_{\mu=0}^{m}
\sum\limits_{q=M_{\mu}+1}^{M_{\mu+1}} \mathcal{P}_{q}(s,t)f\neq  \sum\limits_{\mu=0}^{m}\sum\limits_{k=0}^{\psi(\mu)}\sum\limits_{j=M_{k\mu}+1}^{M_{k\mu+1}} \mathcal{P}_{k_{j}}(s,t)f ,\;f\in   \mathfrak{H}.
$$
However, we can produce      a subsequence $\{\xi_{l}\}_{1}^{\infty}\subset \mathbb{N}$ so that   that the left-hand side and the right-hand side became  equal if $m=\xi_{l},\,l\in \mathbb{N}.$
Note that in accordance with \eqref{18}, we have
\begin{equation*}
B\in \mathfrak{S}_{1}\Rightarrow  |\mu_{n}(B) | =o\left(n^{-1}\right),
\end{equation*}
therefore
$
n(r,B)\leq \varepsilon r,
$
for   sufficiently large values $r,$ where $\varepsilon>0$ is an arbitrary small positive number.
Consider a  counting function  corresponding to the sequence $ \{R_{\mu}\}_{0}^{\infty},$  it is clear that
$
 n(r,R_{\mu})=a^{-1}r  +o(r).
$
Consider the difference
\begin{equation}\label{26}
n(r,R_{\mu})-n(r,B)\geq \left\{a^{-1} -\varepsilon \right\}r  +o(r),
\end{equation}
it is clear that  $n(r,R_{\mu})-n(r,B) \rightarrow\infty,\,r\rightarrow \infty.$ Thus,  we can extract a subsequence  of the natural numbers $\{\xi_{l}\}_{1}^{\infty}\subset \mathbb{N}$  so that the sequence
$
\left\{n(R_{\xi_{l}+1},R_{\mu})-n(R_{\xi_{l}+1},B)\right\}_{1}^{\infty}
$
is monotonically increasing. Therefore   each ring
$\left\{ z\in \mathbb{C}:\;R_{\xi_{l}}<|z|<R_{\xi_{l}+1}\right\}$   does not contain characteristic numbers of the operator $B.$   Now if we consider   possible arrangements of the radii $\tilde{R}_{k\mu}$ then the following fact becomes    clear
$$
M_{\xi_{l}+1} = \sum\limits_{k=0}^{\psi(\xi_{l})}M_{k\xi_{l}+1},\;l=1,2,...\,.
$$
Therefore
$$
\frac{1}{2\pi i}\!\!\!\oint\limits_{\vartheta( R_{\xi_{l}+1}) }  e^{- \lambda^{s}  t}B\left(I-\lambda B\right)^{-1}f d\lambda=-\sum\limits_{\mu=0}^{\xi_{l}}
\sum\limits_{q=M_{\mu}+1}^{M_{\mu+1}} \mathcal{P}_{q}(s,t)f=
$$
$$
=-\sum\limits_{\mu=0}^{\xi_{l}}\sum\limits_{k=0}^{\psi(\mu)}\sum\limits_{j=M_{k\mu}+1}^{M_{k\mu+1}} \mathcal{P}_{k_{j}}(s,t)f, \;f\in   \mathfrak{H},\;l\in  \mathbb{N}.
$$
Apparently, if we prove the fact
$$
\sum\limits_{\mu=0}^{\infty}\sum\limits_{k=0}^{\psi(\mu)}\left\|\sum\limits_{j=M_{k\mu}+1}^{M_{k\mu+1}} \mathcal{P}_{k_{j}}(s,t)f\right\|_{\mathfrak{H}}<\infty,
$$
then we obtain
\begin{equation}\label{27}
\frac{1}{2\pi i}\int\limits_{\vartheta  }  e^{- \lambda^{s}  t}R(\lambda)f d\lambda= \sum\limits_{\mu=0}^{\infty}
\sum\limits_{q=M_{\mu}+1}^{M_{\mu+1}} \mathcal{P}_{q}(s,t)f,
\end{equation}
where the  integration direction is chosen so that the inside of the domain containing the real axis   appears at the right-hand  side while the point is going along the contour,
moreover
$$
\sum\limits_{l=1}^{\infty}\left\|\sum\limits_{\mu=\xi_{l}}^{\xi_{l+1}}
\sum\limits_{q=M_{\mu}+1}^{M_{\mu+1}} \mathcal{P}_{q}(s,t)f\right\|_{\mathfrak{H}}<\infty.
$$
For this purpose
let us estimate $\psi(\mu),$
observe the following relation
$$
\psi(\mu)\leq \mathrm{card}\{\nu\in \mathbb{N}_{0}: |z_{N_{\nu}}|\leq R_{\mu+1}\}+1 <\mathrm{card}\{j\in \mathbb{N} : |\lambda_{j}|\leq R_{\mu+1}\}+1 =
$$
$$
 =n(R_{\mu+1},B)+1 <n(R_{\mu+1},R_{\mu})+1 =\mu+2,
$$
it holds for a sufficiently large value $\mu,$  in accordance with \eqref{26}. Therefore, if we prove that     the following series  is convergent, i.e.
$$
S  \leq\sum\limits_{\mu=0}^{\infty}\sum\limits_{k=0}^{ \mu+2 }\left\|\sum\limits_{j=M_{k\mu}+1}^{M_{k\mu+1}} \mathcal{P}_{k_{j}}(s,t)f\right\|_{\mathfrak{H}} <\infty,
$$
then we obtain the desired result.
Since the inner sum contains the projectors corresponding to the operator $B_{k}$ then we can apply Lidskii V.B. method \cite{firstab_lit:1Lidskii} modified by virtue of Lemma \ref{L9}. Let us estimate the inner sum,  we have
$$
  \left\|\sum\limits_{j=M_{k\mu}+1}^{M_{k\mu+1}} \mathcal{P}_{k_{j}}(s,t)f\right\|_{\mathfrak{H}}=\left\| \,\int\limits_{\vartheta_{k\mu}}  e^{- \lambda^{s}  t}B_{k}\left(I-\lambda B_{k}\right)^{-1}f d\lambda \right\|_{\mathfrak{H}}\leq J_{k\mu}+J_{k\mu+1}+J^{+}_{k\mu}+J^{-}_{k\mu},
$$
$$
\;J_{  k\mu  }: =\left\|\,\int\limits_{\tilde{\vartheta}_{k\mu}} e^{- \lambda^{s}  t}B_{k}\left(I-\lambda B_{k}\right)^{-1}f d\lambda\,\right\|_{\mathfrak{H}},\;J^{+}_{k\mu}: =\left\|\,\int\limits_{\vartheta_{k\mu_{+}}} e^{- \lambda^{s}  t}B_{k}\left(I-\lambda B_{k}\right)^{-1}f d\lambda\,\right\|_{\mathfrak{H}},\;
$$
$$
J^{-}_{k\mu}:= \left\|\,\int\limits_{\vartheta_{k\mu_{-}}} e^{- \lambda^{s}  t}B_{k}\left(I-\lambda B_{k}\right)^{-1}f d\lambda\,\right\|_{\mathfrak{H}},
$$
$$
 \tilde{\vartheta}_{ k\mu }:=\{z\in  \mathbb{C} :  |z|=\tilde{R}_{k\mu},\,|\arg z|\leq \theta_{0}  \},\;\vartheta_{k\mu_{\pm}}:=
\{z\in  \mathbb{C} : \tilde{R}_{k\mu}<|z|<\tilde{R}_{k\mu+1},\, \arg z  =\pm\theta_{0} \}.
$$
$$
\vartheta_{k\mu}=\tilde{\vartheta}_{ k\mu }\cup\vartheta_{k\mu_{+}}\cup \tilde{\vartheta}_{ k\mu+1 }\cup\vartheta_{k\mu_{-}}.
$$
Applying \eqref{24},   we get
$$
 J_{  k\mu  }   \leq \,\int\limits_{\tilde{\vartheta}_{ k\mu }}e^{- t \mathrm{Re}\lambda^{s}}\left\|B_{k}\left(I-\lambda B_{k}\right)^{-1}f \right\|_{\mathfrak{H}}  |d \lambda| \leq
 C\|f\|  e^{ C_{0}   h_{k}  (2eR_{\mu}) \ln\mu }\tilde{R}_{k\mu} \int\limits_{-\theta_{0}}^{\theta_{0}} e^{- t \mathrm{Re}\lambda^{s}} d \,\mathrm{arg} \lambda,\,|\lambda|=\tilde{R}_{k\mu}.
$$
Since in accordance with the made assumptions, we have
$
\theta_{0} < \pi/2s
$
then
$$
\mathrm{Re }\lambda^{s}\geq|\lambda|^{s} \cos  \theta_{0}  s > |\lambda|^{s} \cos \left[(\pi/2s-\varepsilon)s\right]= |\lambda|^{s} \sin \varepsilon s ,\,\lambda \in \vartheta_{k \mu },
$$
where   $\varepsilon$ is a sufficiently small positive value.
Using this estimate,  we get
$$
\ln J_{k\mu}\leq \ln C +  C_{0}    h_{k} (2e R_{\mu})\ln\mu - \tilde{R}_{k\mu}^{s} t \sin \varepsilon s \leq \ln C +  C_{0}    h_{k} (2e R_{\mu})\ln\mu - R_{\mu-1}^{s} t \sin \varepsilon s=
$$
$$
=  \ln C +  R_{\mu-1}^{s}\left\{C_{0}   R_{\mu-1}^{-s} h_{k} (2eR_{\mu})\ln\mu -  t \sin \varepsilon s\right\},\,k\in \mathbb{N}_{0}.
$$
Using \eqref{25}, we get $R_{\mu-1}^{-s}  h_{k} (2eR_{\mu})\ln\mu\rightarrow0,\;\mu\rightarrow\infty,$
uniformly with respect to $k\in \mathbb{N}_{0}.$
Therefore
$$
\ln J_{k\mu}\leq C e^{-C_{1}R_{\mu}^{s}}.
$$
To estimate other terms, we are rather satisfied with  the estimate represented in Lemma 4 (Lidskii V.B.) \cite{firstab_lit:1Lidskii}, what gives us the following relation
$$
\|(I-\lambda B_{k})^{-1}\|\leq \frac{1}{\sin(|\psi|-\theta)},\, \lambda\in\{z\in \mathbb{C}:\,  \arg  z =\psi\},\,\theta<|\psi|<\pi/2.
$$
 Absolutely analogously  to the reasonings represented in Lemma 7   \cite{firstab_lit:1Lidskii}, we get
$$
 J^{\pm}_{k\mu} \leq C\|f\|_{\mathfrak{H}}  \!\!\! \int\limits_{ \tilde{R}_{k\mu} }^{\tilde{R}_{k\mu+1} }  e^{- t \mathrm{Re }\lambda^{s}}   |d   \lambda|\leq C      \int\limits^{\infty}_{ R_{\mu-1} } e^{-t r ^{s} \sin     s \varepsilon}    d   r=s^{-1}\int\limits^{\infty}_{ R^{s}_{\mu-1} } e^{-t \varphi \sin     s \varepsilon}  \varphi^{\frac{1}{s}-1}   d   \varphi\leq \frac{e^{-t R^{s}_{\mu-1}  \sin     s \varepsilon}}{t \sin     s \varepsilon}.
 $$
Thus, we have come to the relation
$$
\left\|\sum\limits_{j=M_{k\mu}+1}^{M_{k\mu+1}} \mathcal{P}_{k_{j}}(s,t)f\right\|_{\mathfrak{H}}\leq C e^{-C_{1}R_{\mu}^{s}}+C e^{-C_{1}R_{\mu+1}^{s}}+2\frac{e^{-t R^{s}_{\mu-1}  \sin     s \varepsilon}}{t \sin     s \varepsilon}.
$$
Combining the obtained estimates, taking into account the fact  $R_{\mu}^{s}\geq C \mu^{s},$ we get

$$
S  \leq\sum\limits_{\mu=0}^{\infty}\sum\limits_{k=0}^{ \mu+2 }\left\|\sum\limits_{j=M_{k\mu}+1}^{M_{k\mu+1}} \mathcal{P}_{k_{j}}(s,t)f\right\|_{\mathfrak{H}} \leq C \sum\limits_{\mu=0}^{\infty}\sum\limits_{k=0}^{ \mu+2 }e^{-C_{2} \mu^{s} }
=\sum\limits_{\mu=0}^{\infty}(\mu+3)e^{-C_{2} \mu^{s} }<\infty.
$$
The convergence of the last series can be verified easily due to the integral test of convergence. Therefore, relation \eqref{27} holds.   Applying   Lemma \ref{L3a}, we complete the proof.

\end{proof}

\subsection{Supplementary remarks and mathematical applications }

\noindent{\bf Evolution equations}\\

The obtained  results admit  the following application.  By virtue of  arbitrariness in choosing  the summation order, we can find   a solution  analytically  for the Cauchy problem for   the fractional evolution equation containing in the second term an operator belonging to a sufficiently wide class, see \eqref{28}. Some  examples  are represented in the paper  \cite{firstab_lit:2kukushkin2022}, wherein  fractional integro-differential operators generated by a $C_{0}$   semigroup of contractions such   as the Riesz potential,  the Riemann-Liouville fractional differential operator, the Kipriyanov operator,  the difference operator are considered.
  General approach, realized  in the paper \cite{kukushkin2021a}, allows us to consider  a special transform of   m-accretive operator. This  approach seems to be  extremely relevant due to   the fact that  the  class of  m-accretive operators  contains the  infinitesimal generator of  a $C_{0}$   semigroup of contractions.   We should remark    that a fractional differential operator  of the real order can be expressed in terms of the  infinitesimal generator of the corresponding semigroup, see \cite{kukushkin2021a}. Below, we represent an abstract scheme of possible applications.

 Consider      element-functions of the Hilbert space  $u:\mathbb{R}_{+}\rightarrow \mathfrak{H},\,u:=u(t),\,t\geq0$    assuming  that if $u$ belongs to $\mathfrak{H}$    then the fact  holds for all values of the variable $t.$ Notice that under such   assumptions all standard topological  properties as completeness, compactness   etc.  remain correctly defined.  We understand such operations as differentiation and integration in the generalized sense that is caused by the topology of the Hilbert space $\mathfrak{H}.$  The derivative is understood as a    limit
$$
  \frac{u(t+\Delta t)-u(t)}{\Delta t}\stackrel{\mathfrak{H}}{ \longrightarrow}\frac{du}{dt} ,\,\Delta t\rightarrow 0.
$$
Let $t\in  J:=[a,b],\,0< a <b<\infty.$ The following integral is understood in the Riemann  sense as a limit of partial sums
\begin{equation*}
\sum\limits_{i=0}^{n}u(\xi_{i})\Delta t_{i}  \stackrel{\mathfrak{H}}{ \longrightarrow}  \int\limits_{ J }u(t)dt,\,\zeta\rightarrow 0,
\end{equation*}
where $(a=t_{0}<t_{1}<...<t_{n}=b)$ is an arbitrary splitting of the segment $ J ,\;\zeta:=\max\limits_{i}(t_{i+1}-t_{i}),\;\xi_{i}$ is an arbitrary point belonging to $[t_{i},t_{i+1}].$
The sufficient condition of the last integral existence is a continuous property (see\cite[p.248]{firstab_lit:Krasnoselskii M.A.}), i.e.
$
u(t)\stackrel{\mathfrak{H}}{ \longrightarrow}u(t_{0}),\,t\rightarrow t_{0},\;\forall t_{0}\in  J.
$
The improper integral is understood as a limit
\begin{equation*}
 \int\limits_{a}^{b}u(t)dt\stackrel{\mathfrak{H}}{ \longrightarrow} \int\limits_{a}^{c}u(t)dt,\,b\rightarrow c,\,c\in  [0,\infty].
\end{equation*}
Consider a fractional integral in the Riemann-Liouvile sense (see \cite{firstab_lit:samko1987})
$$
I^{\alpha}_{-}f(t)=\frac{1}{\Gamma(\alpha)}\int\limits_{0}^{\infty}f(t+x)x^{\alpha-1}dx,\,\alpha\geq0.
$$
Combining the generalized integro-differential    operations, we can consider a  fractional differential operator
 in the Riemann-Liouvile sense, in the formal form, we have
$$
   \mathfrak{D}^{\alpha}_{\!-}f(t):= \frac{(-1)^{n}}{\Gamma(n-\alpha)}\frac{d^{n}}{d t^{n}}\int\limits_{0}^{\infty}f(t+x)x^{n- \alpha-1}dx,\; \alpha\geq0,\,n=[\alpha]+1.
$$
Thus, we can write
$$
\mathfrak{D}^{\alpha}_{\!-}f(t)= (-1)^{n} \frac{d^{n}}{d t^{n}}\left\{I^{n-\alpha}_{-}f(t)\right\}.
$$
Here, we should remark that
$$
 \mathfrak{D}^{n}_{-}f(t)=  (-1)^{n}\frac{d^{n}\!f }{dt^{n}},\; I^{n}_{-}f(t)= \! \int\limits_{t}^{\infty}\! dx_{1}\int\limits_{x_{1}}^{\infty}dx_{2}...\!\!\!\int\limits_{x_{n-1}}^{\infty}\!\!\!f(x_{n})dx_{n},\;n\in \mathbb{N}.
$$
In accordance with the accepted unified form of notation, that stresses the inverse nature of the operators, we can write
$$
\mathfrak{D}^{-\alpha}_{-}f(t)=I^{\alpha}_{-}f(t),\;\alpha \in \mathbb{R}.
$$
Throughout this paragraph, we consider an operator with discrete spectrum   $W$ satisfying the conditions  $\Theta(W)\in \mathfrak{L}_{0}(\theta),\,\theta<\pi/2.$ Applying  the reasonings  represented in   Lemma \ref{L5}, it is not hard to prove that  the operator $W$ is m-accretive. Therefore,      using relations (3.41), (3.53) \cite{firstab_lit:kato1980}, we can  define  fractional powers of the operator $W$ as follows
$$
W^{-\beta}\!f=    \frac{1}{2\pi i}\int\limits_{ \vartheta } \! \lambda^{-\beta} \left(W-\lambda I\right)^{-1}\!  f d\lambda,\;W^{ \beta}\!f=    \frac{1}{2\pi i}\int\limits_{ \vartheta } \! \lambda^{ \beta-1} \left(W-\lambda I\right)^{-1}\!W f d\lambda,
$$
$$
\,f\in \mathrm{D}(W),\,\beta\in (0,1),
$$
where the contour $\vartheta$  is defined in \eqref{22b}.
Let us study   a Cauchy problem
\begin{equation}\label{28}
   \mathfrak{D}^{ \alpha}_{-}  u(t)=  W  u(t),\;\alpha>0,
\end{equation}
with  the initial conditions
$$
\lim\limits_{t\rightarrow 0} \mathfrak{D}^{ k+\alpha-n}_{-}  u(t)= W^{\alpha_{k}}\! h,\;    h\in  \mathrm{D} (W),
$$
$$
 \alpha_{k}=  (\{\alpha\}+k-1)/ \alpha,\;
k=\left\{ \begin{aligned}
  0,1,...,n-1,\;\{\alpha\}\neq 0  \\
  1,2,...,n-1,\;\{\alpha\}= 0   \\
\end{aligned}
 \right.  .
 $$
 Thus, in the case corresponding to  $\alpha=1$    the   Cauchy problem can be rewritten in the classical form
\begin{equation}\label{29}
   \frac{d u}{dt }   =  W  u,\;\;\lim\limits_{t\rightarrow 0} u(t)    =h \in  \mathrm{D} (W).
\end{equation}
This case  was studied  by Lidskii V.B. in the paper  \cite{firstab_lit:1Lidskii},  under the  assumption $B\in \mathfrak{S}_{\sigma},\,\sigma\leq1$ as the most relevant application of the method $(A,\lambda,s).$
In the case $\alpha=2$     the  Cauchy problem can be rewritten in the  form
$$
   \frac{d^{2} u}{dt^{2} }   = -  W   u,\;\lim\limits_{t\rightarrow 0} u(t)    =h,\;\lim\limits_{t\rightarrow 0}  \frac{du}{dt}    =\sqrt{W} h,\; h\in  \mathrm{D} (W).
$$
However, the principal result  obtained in this paper    allows   to consider higher orders  of fractional derivatives independently on the Schatten-von Neumann index  what is reflected in the following theorem.
\begin{teo}\label{T2}
Assume that  the operator $B:=W^{-1}$ satisfies  conditions of Theorem \ref{T1},   $2\theta/\pi<\alpha<1/\sigma\phi,$      then
there exists a solution of the Cauchy problem \eqref{28}
in the form
\begin{equation}\label{29}
u(t)=\sum\limits_{\mu=0}^{\infty}
\sum\limits_{q=M_{\mu}+1}^{M_{\mu+1}} \mathcal{P}_{q}(\alpha^{-1},t)h.
\end{equation}
\end{teo}
\begin{proof}
Consider an element-function
$$
u(t):= \frac{1}{2\pi i}\int\limits_{ \vartheta } e^{- \lambda^{1/\alpha}  t}B\left(I-\lambda B\right)^{-1}h\, d\lambda,\,t>0.
$$
Note that in accordance with  Theorem \ref{T1} relation \eqref{29} holds.  Thus, we should prove   the    fact  that $u(t)$ is a solution of the equation satisfying the initial conditions. In accordance with the preliminary information given  above, we can write
$$
\mathfrak{D}^{ k+\alpha-n}_{-}  u(t)=(-1)^{k} \frac{d^{k}}{d t^{k}}\left\{I^{n-\alpha}_{-}u(t)\right\},\;k=0,1,...,n.
$$
Changing the order of integration, what is based on the statements  of the ordinary calculus, we get
$$
 \Gamma(n-\alpha)\,I^{n-\alpha}_{-}u(t)  = \left(\int\limits_{0}^{\infty}x^{[\alpha]-\alpha}u(t+x)dx, g\right)_{\!\!\mathfrak{H}}=
\frac{1}{2\pi i}\int\limits_{\vartheta }e^{-\lambda^{1/\alpha}   (t+x) } \left(R(\lambda)h ,g\right)_{\mathfrak{H}} d\lambda\int\limits_{0}^{\infty}x^{[\alpha]-\alpha}
dx=
$$
$$
= \frac{1}{2\pi i}\int\limits_{0}^{\infty}x^{[\alpha]-\alpha}e^{-\lambda^{1/\alpha} x}dx\int\limits_{\vartheta }e^{-\lambda^{1/\alpha}   t } \left(R(\lambda)h ,g\right)_{\mathfrak{H}} d\lambda,\;g\in \mathfrak{H},
$$
where $R(\lambda)=B(I-\lambda B)^{-1}.$ Notice that
$$
\int\limits_{0}^{\infty}x^{[\alpha]-\alpha}e^{-\lambda^{1/\alpha} x}dx=\lambda^{1-  n /\alpha }\Gamma(n-\alpha),\,n=[\alpha]+1,
$$
 therefore
$$
I^{n-\alpha}_{-}u(t)=\frac{1 }{2\pi i}\int\limits_{\vartheta }e^{-\lambda^{1/\alpha}   t } \lambda^{1-  n /\alpha } R(\lambda)h   \,d\lambda.
$$
Hence
$$
\mathfrak{D}^{ k+\alpha-n}_{-}  u(t)= \frac{(-1)^{k} }{2\pi i}\frac{d^{k}}{d t^{k}}\int\limits_{\vartheta }e^{-\lambda^{1/\alpha}  t } \lambda^{1- n/\alpha } R(\lambda)h   \,d\lambda = \frac{1 }{2\pi i} \int\limits_{\vartheta }e^{-\lambda^{1/\alpha}   t }  \lambda^{1-  (n-k) /\alpha }  R(\lambda)h   \,d\lambda,
$$
$$
k=0,1,...,n.
$$
The differentiation  under the integral can be easily  substantiated analogously to the proposition related to the ordinary calculus, the detailed reasonings are represented in the proof of Theorem 4 \cite{firstab_lit(frac2023)}.
Taking into account the fact
$
\lambda  R(\lambda)=(I-\lambda B)^{-1}-I,
$
we get
$$
\mathfrak{D}^{ k+\alpha-n}_{-}  u(t)=\frac{1 }{2\pi i} \int\limits_{\vartheta }e^{-\lambda^{1/\alpha}   t }  \lambda^{1-  (n-k) /\alpha }  R(\lambda)h   \,d\lambda=
$$
$$
 =\frac{1 }{2\pi i} \int\limits_{\vartheta }e^{-\lambda^{1/\alpha}   t }  \lambda^{ -  (n-k) /\alpha }  (I-\lambda B)^{-1}h   \,d\lambda +\frac{h }{2\pi i} \int\limits_{\vartheta }e^{-\lambda^{1/\alpha}   t }  \lambda^{ -  (n-k) /\alpha }  d\lambda.
$$
Note that the second integral equals  to zero due to the analytic property of the subintegral function  in the domain
$$
 \left\{z\in \mathbb{C}: |z|=r_{0}, |z|=R, |\arg z|\leq \theta_{0}\right\}\cup\left\{z\in \mathbb{C}: r_{0}<|z|<R,  |\arg z|=\theta_{0}\right\},
$$
where $R>0$ is an arbitrary large number, and the fact
$$
 \int\limits_{\vartheta_{\! R} }e^{-\lambda^{1/\alpha}   t }  \lambda^{ -  (n-k) /\alpha }  d\lambda\rightarrow 0,\,R\rightarrow\infty,\;\vartheta_{\! R}:=\left\{z\in \mathbb{C}: |z|=R,\,|\arg z|\leq \theta_{0}\right\}.
$$
Therefore
$$
\mathfrak{D}^{ k+\alpha-n}_{-}  u(t)=\frac{1 }{2\pi i} \int\limits_{\vartheta }e^{-\lambda^{1/\alpha}   t }  \lambda^{ -  (n-k) /\alpha }  (I-\lambda B)^{-1}h   \,d\lambda=\frac{1 }{2\pi i} \int\limits_{\vartheta }e^{-\lambda^{1/\alpha}   t }  \lambda^{ -  (n-k) /\alpha }  R(\lambda) Wh   \,d\lambda =
$$
\begin{equation}\label{29a}
=\frac{1 }{2\pi i} \int\limits_{\vartheta }e^{-\lambda^{1/\alpha}   t }  \lambda^{ -  (n-k) /\alpha } W R(\lambda) h   \,d\lambda,
\end{equation}
from what follows the fact  $\mathfrak{D}^{  \alpha }_{-}  u(t)=W u(t),$ if we put $k=n.$ Assume that   $k$ corresponds to the initial conditions.
 Observe  that  the    commutative property
 $
 B(I-\lambda B)^{-1} =(I-\lambda B)^{-1}B,
 $
  see Problem 5.4 \cite[p.36]{firstab_lit:kato1980},  leads   to   the equality $
  (W-\lambda I)^{-1}= R(\lambda).
$
It follows that  that the  improper  integral \eqref{29a} is uniformly convergent with respect to the parameter $t$ since we have
$$
 \|R(\lambda)\|\leq C|\lambda|^{-1},  \, \lambda\in\{z\in \mathbb{C}:\,  \arg  z =\psi\},\,\theta<|\psi|<\pi/2,
 $$
 see the proof of Lemma \ref{L5}. Therefore,    passing to the limit under the integral, we get
$$
\mathfrak{D}^{ \alpha-n}_{-}  u(t)\rightarrow \frac{1 }{2\pi i} \int\limits_{\vartheta }   \lambda^{ 1-   n   /\alpha }  R(\lambda)  h   \,d\lambda =W^{1-   n  /\alpha}h,\,t\rightarrow 0,\,\{\alpha\}\neq0,
$$
$$
\mathfrak{D}^{ k+\alpha-n}_{-}  u(t)\rightarrow \frac{1 }{2\pi i} \int\limits_{\vartheta }   \lambda^{ -  (n-k) /\alpha }  R(\lambda) Wh   \,d\lambda =W^{1-  (n-k) /\alpha}h,\,t\rightarrow 0,
$$
$$
k=1,2,...,n-1,
$$
where the last  formula is true except for the case $\{\alpha\}=0,\,k=1$    which is covered by Lemma \ref{L5}.   The proof is complete.
\end{proof}

The main advantage of the last theorem is  that there are not   any  restrictions upon  the highest  order of the fractional derivative   or the index of the Schatten-von Neumann  class in the case when the operator $W^{-1}$ has the sequence of the algebraic multiplicities of the lowest growth.
  It   follows from the opportunity to  consider an arbitrary small value of the summation order in accordance with    Theorem \ref{T1}.\\

\newpage

\noindent{\bf   Spectral asymptotics for fractional-differential and pseudo-differential operators }\\

In this paragraph, we study  an  operator with   discrete spectrum $W.$    Consider a   formula obtained by  Markus A.S.,  Matsaev V.I. (2.1) \cite{firstab_lit:Markus Matsaev} connecting spectral asymptotic  of the operator and its real component.   Assume that    condition
\begin{equation}\label{11a}
|(\mathfrak{Im}W f,f)_{\mathfrak{H}}|\leq \frac{b}{2} \big\|\sqrt{\mathfrak{Re} W}f\big\|^{2q}_{\mathfrak{H}} \|f\|^{2-2q}_{\mathfrak{H}},
 \,f\in \mathrm{D} \big(\sqrt{\mathfrak{Re} W} \big),\;q\in [0,1),
\end{equation}
holds then in accordance with  Theorem 2.1  \cite{firstab_lit:Markus Matsaev}    for an arbitrary $\delta>0$  and $r>r_{1}$ the following relation holds
\begin{equation}\label{30}
n(r,W)-n(r,\mathfrak{Re}W)\leq K\left\{n\left(r+b[1+\delta]r^{q} ,\mathfrak{Re}W\right)-n\left(r-b[1+\delta]r^{q} ,\mathfrak{Re}W\right)    \right\},
\end{equation}
where $K$ is a constant that depends on $\delta$ only, the constant  $r_{1}$ depends of $\delta,q,b.$ Here, we should remark that  instead of  condition  \eqref{11a} ($2^{\,0}$ in accordance with the terminology used in  \cite{firstab_lit:Markus Matsaev}),   we can impose  more general condition   $3^{\,0} \cite{firstab_lit:Markus Matsaev}$ that  guarantees \eqref{30}.\\

Consider a spectral  asymptotics that is inherent to  a wide class  of fractional-differential and pseudo-differential operators, see \cite{firstab_lit:Rosenblum}
$$
n(r,\mathfrak{Re}W)=\gamma_{0} r^{\xi}+o(r^{\mu}),\,r\rightarrow\infty,\,0\leq\mu<\xi.
$$
Taking into account the fact obtained due  the Taylor   formula
$$
(r+cr^{q})^{\xi}-(r-cr^{q})^{\xi}=r^{\xi}\left\{(1+cr^{q-1})^{\xi}-(1-cr^{q-1})^{\xi}\right\}=  2\xi cr^{\xi+q-1}+O(r^{\xi+2q-2}),
$$
where $c=b[1+\delta],$ substituting, we get
 $$
 n(r,W)-n(r,\mathfrak{Re}W)\leq K\left\{  2\xi cr^{\xi+q-1}+o(r^{\mu})+O(r^{\xi+2q-2})  \right\}.
 $$
 Therefore $n(r,W)-n(r,\mathfrak{Re}W)=o(r^{\mu}),\,\mu> \xi+q-1.$
 The given above reasonings lead  us to the   implication
\begin{equation}\label{31}
n(r,\mathfrak{Re}W)= \gamma_{0}r^{\xi}+o(r^{\mu}),\,\Rightarrow n(r,W)= \gamma_{0}r^{\xi}+o(r^{\mu}),\;\mu> \xi+q-1 .
\end{equation}
 However, we have an interest in a more explicit formula      studied in  \cite{firstab_lit:Duistermaat},\cite{firstab_lit:Rosenblum}
\begin{equation}\label{32}
n(r,\mathfrak{Re}W)=   \sum\limits_{j=0}^{l} \gamma_{j}r^{\frac{v-j}{m}} + o(r^{\frac{v-l}{m}}),\,0\leq l\leq v,
\end{equation}
 where $v\in \mathbb{N}$ is the dimension of the Euclidian space    $m\in \mathbb{N}$ is the order of the operator (derivative),       $\gamma_{j}$ are  constants, $\gamma_{0}>0.$   Observe that by virtue of  the assumption  $l/m<1-q$ we are able   to apply   the  scheme of reasonings used to obtain \eqref{31},
analogously to the  above, we get
\begin{equation}\label{33}
n(r,\mathfrak{Re}W)=   \sum\limits_{j=0}^{l} \gamma_{j}r^{\frac{v-j}{m}} + o(r^{\frac{v-l}{m}}),\Rightarrow n(r,W)= \sum\limits_{j=0}^{l} \gamma_{j}r^{\frac{v-j}{m}} + o(r^{\frac{v-l}{m}}).
\end{equation}
This fact is noticed  in   Corollary 2.3 \cite{firstab_lit:Markus Matsaev} for two-termed asymptotics.
On the other hand,  implementing the scheme of reasonings   represented in  \cite[p.98]{firstab_lit:2Agranovich1990} (6.1.15),   we obtain
\begin{equation}\label{34}
n(r,W)= \sum\limits_{j=0}^{l} \gamma_{j}r^{\frac{v-j}{m}} + o(r^{\frac{v-l}{m}}),\Rightarrow  |\mu_{n}(W)| = \sum\limits_{j=0}^{l}\tilde{\gamma}_{j}n^{\frac{m-j}{v}}+  o(n^{\frac{m-l}{v}}),
\end{equation}
where $\tilde{\gamma}_{j}$ are   constants $\tilde{\gamma}_{0}>0.$
Now, if  we assume that $v=l,$ then
$
 |\mu_{n+1}(W)|- |\mu_{n}(W)|>0
$
for sufficiently large numbers $n\in \mathbb{N}.$ It follows that the operator $W^{-1}$ has the sequence of the  algebraic multiplicities of the lowest growth.
Summarizing the given above information, we conclude that the case $l=v,\,v/m<1-q$ falls within the scope of the developed method. Further, we illustrate the idea through  the low values of $l$ corresponding to   well-known operators.\\

 We have a particular interest in the case when condition \eqref{11a} is satisfied and the following  relation holds
\begin{equation}\label{35}
n(r, \mathfrak{Re} W)=\gamma_{0}r^{\xi}+O( \ln r),\, 0<\xi<1-q,
\end{equation}
  then in accordance with Corollary 2.2 \cite{firstab_lit:Markus Matsaev}, we get
\begin{equation*}\label{36a}
n(r,W)=\gamma_{0}r^{\xi}+O( \ln r).
\end{equation*}
By direct calculations, we obtain
$$
n^{a/\xi}(r)=r^{a}\left\{\gamma_{0}+r^{-\xi}O(\ln r) \right\}^{a/\xi}  =\gamma^{a/\xi}_{0} r^{a}+r^{a-\xi}O(\ln r),\,a>0,
$$
where $n(r):=n(r,W),$
therefore
$$
\frac{r^{1-\xi}O(\ln r)}{n^{a/\xi}(r)}\rightarrow 0,\,r\rightarrow\infty,\,a>1-\xi.
$$
Thus, we obtain
\begin{equation}\label{36}
r=\gamma^{-1/\xi}_{0}n^{1/\xi}(r)+r^{1-\xi}O(\ln r);\,r=\gamma^{-1/\xi}_{0}n^{1/\xi}+ o(n^{a/\xi}),\;n\in \mathbb{N}.
\end{equation}
Now assume that $\xi=v/m,\, 1<v<m,$ then for the sake of certainty,  we can choose $a$    satisfying the   condition
$$
a/\xi=am/v=(m-1)/v>m/v-1=(1-\xi)/\xi,
$$
and rewrite relation \eqref{36} in the form
$$
|\mu_{n}(W)|=\gamma^{-1/\xi}_{0}n^{m/v}+ o(n^{(m-1)/v}).
$$
It follows that for sufficiently large values $j\in \mathbb{N},$ we have
$$
|\mu_{(j+1)^{v}}(W)|-|\mu_{  j ^{v}}(W)|>0,
$$
hence the principal indexes admit the following estimate $p_{j}\leq C j^{v}.$   Consider a function $\Lambda(r):=n(r,W)-n(r-\delta,W),$ where $\delta>0$ is an arbitrary small positive fixed number. It is clear that
$$
\Lambda(r)\leq C r^{v/m-1}+O(\ln r).
$$
Substituting \eqref{36}, taking into account the fact $v<m,$ we obtain
$$
  \Delta_{j}\leq \Lambda(\mu_{p_{j}})\leq C_{1}p_{j}^{1-m/v}+C_{2}\ln C_{3}p_{j}\leq C \ln j,\,j\in \mathbb{N},
$$
from what follows that the operator has the sequence of the algebraic multiplicities of  the lowest growth. It is remarkable that conditions \eqref{33},\eqref{35} allow    to establish the fact that the summation order   is an arbitrary small positive  value  for a sufficiently large operator  class. Some concrete  examples of operators satisfying the conditions     can be found in \cite{firstab_lit:Rosenblum}.\\

  In order to represent a concrete example, consider the one-dimensional Schrodinger operator \cite[p.194]{firstab_lit:Rosenblum}
$$
L=-d^{2}/dx^{2}+\varrho(x),\;\varrho(x)\geq a>0,\;\varrho\in CAP^{\infty}\!(\mathbb{R}).
$$
 The following asymptotical formula was obtained in the paper \cite{firstab_lit:Savin1988}, see also \cite[p.175]{firstab_lit:Rosenblum}
$$
n(r,L)=\sqrt{r}+r^{-\frac{1}{2}}\sum\limits_{k=0}^{\infty}d_{k}r^{-k}.
$$
Consider an operator with a discrete spectrum    $W$ satisfying   the following conditions
\begin{equation}\label{37}
\mathfrak{Re}W=L,\;|\left( \mathfrak{Im}Wf,f \right)_{L_{2}(\mathbb{R})}|\leq 2^{\,-\frac{1}{2} } \left(Lf,f \right)^{q}_{L_{2}(\mathbb{R})},\,f\in C_{0}^{\infty}(\mathbb{R}),\,\|f\|_{L_{2}(\mathbb{R})}=1,
\end{equation}
$$
0<q<1/2.
$$
Applying Theorem 2.1 \cite{firstab_lit:Markus Matsaev}, i.e. implication  \eqref{33},  then using  implication  \eqref{34},     we obtain
$$
|\mu_{n}(W)|=n^{2} +o(n),\,n\in \mathbb{N},
$$
hence
$
|\mu_{n+1}(W)|-|\mu_{n}(W)|>0
$
for sufficiently large numbers $n\in \mathbb{N}.$
It follows that the operator $B:=W^{-1}$ has the sequence of the  algebraic multiplicities of the lowest growth. The sectorial condition is satisfied  by virtue of relation  \eqref{37}, i.e $\Theta(B)\subset\mathfrak{L}_{0}(\theta),\,\theta<\pi/4.$
Applying Lemma 1 \cite{firstab_lit(Math2024)} to the  compact sectorial operator $B$   we get
\begin{equation*}
s_{2n-1}(B)\leq \sqrt{2}\,\sec  \theta \cdot \mu_{n} (\mathfrak{Re}  B),\;\;s_{2n}(B)\leq \sqrt{2}\,\sec  \theta \cdot  \mu_{n} (\mathfrak{Re}  B), \; n\in \mathbb{N}.
\end{equation*}
Using the  properties of the operator $L,$ condition \eqref{37}, we can prove that the operator $W$ satisfies  conditions $\mathrm{H}1,\mathrm{H}2$ \cite{kukushkin2021a}.  In accordance with the results  \cite{firstab_lit(arXiv non-self)kukushkin2018} the conditions  $\mathrm{H}1,\mathrm{H}2$  guarantee
$$
\mu_{n}^{-1}(L)\asymp  \mu_{n}(\mathfrak{Re}B).
$$
Analogously to  \eqref{34}, we get
$$
 \mu_{n}(L) =n^{2} +o(n),\,n\in \mathbb{N}.
$$
Combining these relations, we obtain the fact  $B\in \mathfrak{S}_{1}.$ Thus, the operator $B$ satisfies conditions of Theorem \ref{T1} and therefore operator $W$ satisfies conditions of Theorem \ref{T2}. This example illustrates efficiency of the obtained results in studying evolution equations containing a non-selfadjoint term. \\

  Consider  an operator class   generated by a $C_{0}$ semigroup  of contractions \cite{kukushkin2021a}.   We can refer Theorem 5  \cite{kukushkin2021a} in accordance with which the operators belonging  to the class  satisfy conditions H1,H2 \cite{kukushkin2021a} and therefore have convenient, from the created theory point of view, properties such as a sectorial property, compactness of the resolvent, belonging to the Schatten-von Neumann  class, etc. Consider the infinitesimal generator $A$ of a $C_{0}$ semigroup of contractions, we can form the infinitesimal  generator transform
$$
 Z_{\alpha}(A):= A^{\ast}GA+FA^{\alpha},\,\alpha\in [0,1),
$$
where the symbols  $G,F$  denote  operators acting in $\mathfrak{H}.$   Taking into account Corollary 3.6 \cite[p.11]{Pasy},  Theorem 5 \cite{kukushkin2021a}, we conclude that if    $ A^{-1}$ is compact,  $F,G\in \mathcal{B}(\mathfrak{H}),$  $G$ is    strictly accretive,
then   $Z_{\alpha} (A)$ satisfies  conditions  H1-H2 \cite{kukushkin2021a}.  Note that   Theorem 5 \cite{kukushkin2021a} gives us a tool to describe spectral properties  of the  transform $Z_{\alpha} (A),$  in particular, we can establish the  index of  the  Schatten-von Neumann class  applying   Theorem 3 \cite{kukushkin2021a}. Apparently, having known the index of the Schatten-von Neumann  class, we can proceed to the next step   applying  results of this paragraph in order to verify fulfilment of Theorem 1 conditions.

\section{Conclusions}

   In the paper, we  have  shown that the summation order in the Abel-Lidskii sense  can be decreased   to  an arbitrary small positive value in the case corresponding to a sectorial operator belonging to the trace class under conditions imposed upon the growth of the algebraic multiplicities.  In addition, we  produce  a number of fundamental  propositions in the framework of functional analysis  which may represent the interest themselves. The lemma   on estimation of the characteristic determinant and the lemma on discharging of the spectrum  create  a prerequisite  for further study. Thus,  the invented  technique forms a base for extension of the obtained  results to an arbitrary Schatten-von Neumann class. Moreover, based upon the general scheme of reasonings and  studying more detailed the issue on  the   brackets arrangement     in the series, we can  construct a qualitative theory  describing the peculiarities of the  series summation in the Abel-Lidkii sense. The application part appeals to the existence  theorem for the abstract Cauchy problem for the fractional evolution equation that covers  many  concrete problems in the theory of differential equations. The  results can be clearly illustrated on     the operators having  sufficiently  slow growth of algebraic multiplicities. In its own turn,   the latter property can be  clearly expressed through the scale of spectral asymptotics for selfadjoint operators which is studied in detail in the final paragraph. In author's opinion, the significant achievement of the paper is    that  the  operator class  generated by   strictly continuous semigroups of contractions  can be studied due to the obtained methods. The latter  includes  many well-known integro-differential operators such  as  the linear combination of the  differential operator and the Kipriyanov fractional-differential operator, the linear combination of the  differential operator and the  Riesz potential, the perturbation of the  difference operator.
 The author reasonably believes  that the represented   statements  are principally novel in the framework of the abstract   spectral theory while the obtained conclusions  admit  relevant  applications.


\newpage







 \begin{thebibliography}{2}


\bibitem{firstab_lit:Agranovich1976} {\sc  Agranovich M.S.} Summability of series in root vectors of non-selfadjoint elliptic operators. \textit{ Funct. Anal. Appl.}, \textbf{10}  (1976),    165-174





\bibitem{firstab_lit:2Agranovich1990} {\sc  Agranovich M.S.} Elliptic operators on closed manifolds.
 \textit{Partial differential equations. VI. Elliptic operators on closed manifolds. Encycl. Math. Sci.; translation from Itogi Nauki Tekh., Ser. Sovrem. Probl. Mat., Fundam. Napravleniya},  \textbf{63}  (1990), 5-129.





\bibitem{firstab_lit:2Agranovich1994} {\sc  Agranovich M.S.} On series with respect to root vectors of operators associated with forms having symmetric principal part.
\textit{Functional Analysis and its applications},   \textbf{28} (1994),   151-167.


\bibitem{firstab_lit: 1999 Katsenelenbaum} {\sc  Agranovich M.S.,  Katsenelenbaum B.Z., Sivov A.N., Voitovich N.N.}
Generalized method of eigenoscillations in the diffraction theory. \textit{Zbl 0929.65097 Weinheim: Wiley-VCH.}, 1999.











 \bibitem{firstab_lit:Bazhl} {\sc Bazhlekova  E.} The abstract Cauchy problem for the fractional evolution
equation.
\textit{Fractional Calculus and Applied Analysis},   \textbf{1}, No.3 (1998),  255-270.


\bibitem{firstab_lit:Bazhl1} {\sc Bazhlekova  E.}  The abstract Cauchy problem for the fractional evolution equation. \textit{RANA: reports on
applied and numerical analysis, Technische Universiteit Eindhoven}, \textbf{9814} (1998), 1-17.

\bibitem{firstab_lit:Bazhl2} {\sc Bazhlekova  E.}   Fractional evolution equations in Banach spaces. \textit{Eindhoven: Technische Universiteit
Eindhoven}, 2001; https://doi.org/10.6100/IR549476.




\bibitem{firstab_lit:Baraichev} {\sc Braichev G.G.}  Exact relationships between certain characteristics of growth for complex sequences.\textit{Ufa Mathematical Journal}, \textbf{5},  No.4   (2013), 16-29.















\bibitem{firstab_lit: Courant-Hilbert}{\sc Courant R., Hilbert D.}  Methods of mathematical physics.\textit{Gostekhizdat, Moscow},   1951.




\bibitem{firstab_lit:Duistermaat} {\sc Duistermaat L.J., Guillemin V.} The spectrum of positive elliptic operators
and periodic bicharacterisitids.
\textit{Invent math.}, \textbf{29}, No.1   (1975), 39-79.




\bibitem{firstab_lit:1Gohberg1965} {\sc   Gohberg I.C., Krein M.G.}  Introduction to the theory of linear non-selfadjoint operators in a Hilbert space.
 \textit{Moscow: Nauka, Fizmatlit},  1965.







\bibitem{firstab_lit:Hardy}
 {\sc Hardy G.H.} Divergent series.
  \textit{Oxford University Press, Ely House, London W.}, 1949.







\bibitem{firstab_lit:kato1980}{\sc Kato T.} Perturbation theory for linear operators. \textit{Springer-Verlag Berlin, Heidelberg, New York}, 1980.




\bibitem{firstab_lit:1967Katsnelson} {\sc Katsnelson V.E.} Convergence and summability of series in the root vectors of certain classes of
non-selfadjoint operators.
\textit{Ph.D thesis, Kharkiv State University, Kharkiv}, 1967. (in Russian)











\bibitem{firstab_lit:1Katsnelson} {\sc Katsnelson V.E.} Conditions under which systems of eigenvectors of some classes
of operators form a basis.
\textit{Funct.Anal. Appl.}, \textbf{1}, No.2   (1967), 122-132.








\bibitem{firstab_lit:Keldysh M.V.}{\sc Keldysh M.V.}  On eigenvalues and eigenfunctions of some classes of non-selfadjoint equations. \textit{Doklady Akademii Nauk SSSR}, 77  (1951), 11-14.



\bibitem{firstab_lit:Koch1909}
{\sc Von Koch H.}  Sur la convergence des determinants infinies. \textit{Rend. Circ Math. Palermo}, \textbf{28}  (1909), 255-266.


















\bibitem{firstab_lit(arXiv non-self)kukushkin2018}   {\sc  Kukushkin M.V.} On One Method of Studying Spectral Properties of Non-selfadjoint Operators. \textit{Abstract and Applied Analysis; Hindawi: London, UK}, \textbf{2020}   (2020);    https://doi.org/10.1155/2020/1461647.


\bibitem{kukushkin2019}{\sc Kukushkin M.V.} Asymptotics of eigenvalues for differential
operators of fractional order.\textit{Fract. Calc. Appl. Anal.} \textbf{22}, No. 3 (2019), 658-681; DOI:10.1515/fca-2019-0037; at https://www.degruyter.com/view/j/fca.


\bibitem{kukushkin2021a}{\sc Kukushkin M.V.} Abstract fractional calculus for m-accretive operators.
 \textit{International Journal of Applied Mathematics.} \textbf{34}, Issue: 1 (2021);  DOI: 10.12732/ijam.v34i1.1.








\bibitem{firstab_lit:1kukushkin2021}   {\sc  Kukushkin M.V.} Natural lacunae method and Schatten-von Neumann classes of the convergence exponent. \textit{Mathematics}, \textbf{10}, (13), 2237   (2022); https://doi.org/10.3390/math10132237.

\bibitem{firstab_lit:2kukushkin2022}   {\sc Kukushkin M.V.} Evolution Equations in Hilbert Spaces via the Lacunae Method. \textit{Fractal Fract.}, \textbf{6}, (5), 229   (2022); https://doi.org/10.3390/fractalfract6050229.

\bibitem{firstab_lit(axi2022)}   {\sc Kukushkin M.V.} Abstract Evolution Equations with an Operator Function in the Second Term.
 \textit{Axioms}, 11, 434  (2022);   https://doi.org/10.3390/axioms 11090434.






\bibitem{firstab_lit(frac2023)}{\sc  Kukushkin M.V.}  Cauchy Problem for an Abstract Evolution Equation of Fractional Order. \textit{Fractal Fract.}, 7, 111 (2023);     https://doi.org/10.3390/fractalfract7020111.


\bibitem{firstab_lit(Math2024)}  {\sc Kukushkin M.V.}, Schatten
Index of the Sectorial Operator via the
Real Component of Its Inverse.   \textit{Mathematics},   \textbf{12}, 540   (2024); https://doi.org/10.3390/math12040540.


\bibitem{firstab_lit:Krasnoselskii M.A.}{\sc Krasnoselskii M.A.,  Zabreiko P.P.,  Pustylnik E.I.,  Sobolevskii P.E.}
  Integral operators in the spaces of summable functions. \textit{ Moscow:  Science,   FIZMATLIT},   1966.


\bibitem{firstab_lit:1959Krein}{\sc  Krein M.G.} Criteria for completeness of a system of root vectors of a dissipative operator.
  \textit{Uspekhi Matematicheskikh Nauk}, \textbf{14}, Issue 3(87) (1959), 145-152.



\bibitem{firstab_lit:1Krein} {\sc  Krein M.G.}  Criteria for
completeness of the system of root vectors of a dissipative operator.
 \textit{Amer. Math. Soc. Transl. Ser., Amer. Math. Soc., Providence, RI}, \textbf{26}, No.2 (1963), 221-229.




\bibitem{firstab_lit:Eb. LevinE} {\sc  Levin  B. Ja.}  Lectures on  Entire Functions.
 \textit{Translations of Mathematical Monographs},     1991.

\bibitem{firstab_lit:Eb. Levin} {\sc  Levin  B. Ja.}  Distribution of Zeros of Entire Functions.
 \textit{Translations of Mathematical Monographs},     1964.







\bibitem{firstab_lit:1958Lidskii} {\sc  Lidskii V.B.} Theorems on the completeness of a system of characteristic and adjoined elements of operators having a discrete spectrum.
\textit{Dokl. Akad. Nauk SSSR }, \textbf{119}, No. 6  (1958), 1088-1091.


\bibitem{firstab_lit:1959Lidskii} {\sc  Lidskii V.B.} Non-selfadjoint operators with a trace.
 \textit{Dokl. Akad. Nauk SSSR },\textbf{125}, No. 3 (1959), 485-487.

\bibitem{firstab_lit:1960Lidskii} {\sc  Lidskii V.B.}   Summation of series over the main vectors of non-selfadjoined operators. \textit{Soviet Math. Dokl.}
  \textbf{132},   No. 2  (1960), 275-278.

\bibitem{firstab_lit:1Lidskii} {\sc  Lidskii V.B.} Summability of series in terms of the principal vectors of non-selfadjoint operators.
 \textit{Tr. Mosk. Mat. Obs.}, \textbf{11}  (1962), 3-35.

\bibitem{firstab_lit:3Lidskii} {\sc  Lidskii V.B.} On the Fourier series expansion on the major vectors  of a non-selfadjoint elliptic operator.
 \textit{Tr. Mosk. Mat. Obs.}, \textbf{57} (99)  (1962), 137-150. (in Russian)

\bibitem{firstab_lit:2Lidskii} {\sc  Lidskii V.B.}  Conditions for completeness of a system of
root subspaces for non-selfadjoint operators with discrete spectra.  \textit{Amer. Math. Soc. Transl. Ser., Amer. Math. Soc., Providence, RI.},
\textbf{34}, No. 2 (1963),  241-281.











\bibitem{firstab_lit:1Markus}  {\sc Markus A.S.}   On   the  basis of root vectors
of a dissipative operator. \textit{Soviet Math. Dokl.},  \textbf{1}   (1960), 599-602.


\bibitem{firstab_lit:2Markus} {\sc Markus A.S.} Expansion in root vectors of a slightly perturbed selfadjoint operator.
\textit{ Soviet Math. Dokl.}, \textbf{3}   (1962), 104-108.

\bibitem{firstab_lit:1966Markus}   {\sc Markus A.S.} Certain criteria for the completeness of a system of root vectors of a linear operator in a Banach space.
\textit{Matematicheskii Sbornik. Novaya Seriya}, \textbf{70} (112), No. 4  (1966),    526-561



\bibitem{firstab_lit:Markus Matsaev} {\sc Markus A.S.,  Matsaev V.I.} Operators generated by sesquilinear forms and their spectral asymptotics.
\textit{ Linear operators and integral equations, Mat. Issled., Stiintsa, Kishinev}, \textbf{61} (1981), 86-103.


\bibitem{firstab_lit:Markus Matsaev1} {\sc Markus A.S.,  Matsaev V.I.} On the convergence of the root vector series for operators close to  selfadjoint ones.
\textit{ Linear operators and integral equations, Mat. Issled., Stiintsa, Kishinev}, \textbf{61}  (1981), 103-129.



\bibitem{firstab_lit:Motovilov} {\sc   Motovilov A. K.,  Shkalikov A. A.}
 Preserving of the unconditional basis property under non-self-adjoint perturbations of self-adjoint operators.
 \textit{Funktsional. Anal. i Prilozhen.}, \textbf{53}, Issue 3  (2019),  45-60.



\bibitem{Pasy} {\sc Pazy A.} Semigroups of Linear Operators and Applications to Partial Differential Equations.
  \textit{Berlin-Heidelberg-New York-Tokyo, Springer-Verlag (Applied Mathematical Sciences V. 44)},  1983.



\bibitem{firstab_lit:Riesz1955} {\sc  Riesz F.,  Nagy B. Sz.} Functional Analysis.
  \textit{Ungar, New York}, 1955.











\bibitem{firstab_lit:Rosenblum} {\sc  Rozenblyum G.V.,  Solomyak M.Z., Shubin M.A.} Spectral theory of differential operators.
 \textit{Results of science and technology. Series Modern problems of mathematics
Fundamental directions}, \textbf{64}    (1989),    5-242.





 \bibitem{firstab_lit:sadovn} {\sc  Sadovnichii V.A.} Theory of operators.
  \textit{Monographs in Contemporary Mathematics, Springer New York, NY}, 1991.

\bibitem{firstab_lit:1959Sakhnovich}{\sc  Sakhnovich L.A.}  A study of the  triangular form  of non-selfadjoint operators.
  \textit{Izvestiya Vysshikh Uchebnykh Zavedenii. Matematika}, No. 4 (1959), 141-149.



\bibitem{firstab_lit:samko1987} {\sc Samko S.G., Kilbas A.A., Marichev O.I.} Fractional Integrals and Derivatives: Theory and Applications.
  \textit{Gordon and Breach Science Publishers: Philadelphia, PA, USA}, 1993.



\bibitem{firstab_lit:Savin1988} Savin A.V. Asymptotic expansion of the density of states for one-dimensional Schrodinger and Dirac operators with almost periodic and random potentials. \textit{Collection of scientific papers. M}, 1988. (in Russian)






\bibitem{firstab_lit:1982 Shkalikov} {\sc  Shkalikov A.A.}  Estimates of meromorphic functions and summability theorems.
 \textit{Pacific J. Math.},   103:2  (1982),  569-582.

\bibitem{firstab_lit:1983 Shkalikov} {\sc  Shkalikov A.A.}  On estimates of meromorphic
functions and summation of series in the root vectors of nonselfadjoint operators.
 \textit{Soviet Math. Dokl.}, \textbf{27}  (1983),  259-263.

\bibitem{firstab_lit:Shkalikov A.} {\sc  Shkalikov A.A.}  Perturbations of selfadjoint and normal operators with a discrete spectrum.
 \textit{Russian Mathematical Surveys}, \textbf{71}, Issue 5 (431) (2016),  113-174.








 \end{thebibliography}
\end{document}